\documentclass[12 pt]{amsart} 
\usepackage{amsfonts} 
\usepackage{amssymb}
\usepackage{amsmath} 
\usepackage{graphicx}
\usepackage{color}

\usepackage[T1]{fontenc}

\theoremstyle{definition} 
\newtheorem{defs}{Definition}[section]

\theoremstyle{plain} 
\newtheorem{prop}[defs]{Proposition}
\newtheorem{lem}[defs]{Lemma}
\newtheorem{ths}[defs]{Theorem} 
\newtheorem{cor}[defs]{Corollary}
\newtheorem*{NNths}{Theorem}
\newtheorem*{NNprop}{Proposition} 
\newtheorem*{NNlem}{Lemma}

\newtheorem*{MT}{Main Theorem}

\theoremstyle{remark} 
\newtheorem{rem}[defs]{Remark}

\newtheorem{Example}[defs]{Example} 
\newtheorem{conv}{Convention}

\title{Degenerations of quadratic differentials on $\mathbb{CP}^1$}
\author{Corentin Boissy}
\address{IRMAR,
Universit\'e Rennes 1,
Campus de Beaulieu,
35042 Rennes cedex, France}
\email{corentin.boissy@univ-rennes1.fr}

\subjclass[2000]{Primary: 32G15. Secondary: 30F30, 57R30}
\keywords{Quadratic differentials,  configurations, \^homo\-logous saddle connections, Siegel-Veech constants}

\begin{document}
                                    

\begin{abstract} 

We describe the connected components  of the complement of a natural ``diagonal'' of real codimension $1$ in  a stratum of quadratic differentials on $\mathbb{CP}^1$.  We establish a natural bijection between the set of these connected components and the set of  generic configurations  that appear on such ``flat spheres''. We also prove that the stratum has only one topological end.
Finally, we elaborate a necessary toolkit destined to evaluation of the Siegel-Veech constants.
\end{abstract}

\maketitle

\setcounter{tocdepth}{1} \tableofcontents

\section{Introduction} 
The article deals with families of flat metric on surfaces of genus zero, where the flat metrics are assumed to have conical singularities,  $\mathbb{Z}/ 2\mathbb{Z}$ linear holonomy and a fixed vertical direction. The moduli space of such metrics is isomorphic to the moduli space of meromorphic quadratic differential on $\mathbb{CP}^1$ with at most simple poles and is naturally stratified by the number of poles and by the orders of zeros of a quadratic differential.

Any stratum is non compact and a neighborhood of its boundary consists of flat surfaces that admit saddle connections of small length. The structure of the neighborhood of the boundary is also related to counting problems in a generic surface of the strata (the ``Siegel-Veech constants, see \cite{EMZ} for the case of Abelian differentials). 

When the length of a saddle connection tends to zero, some other saddle connections might also be forced to shrink.  In the case of an Abelian differential this correspond to homologous saddle connections. In the general case of quadratic differentials, the corresponding collections of saddle connections on a flat surface are said to be \emph{\^homo\-logous}\footnote{The corresponding cycles are in fact homologous on the canonical double cover of $S$, usually denoted as $\widehat{S}$, see section \ref{homolo}.} (pronounced ``hat-homologous''). Configurations associated to collections of \^homo\-logous saddle connections have been described for general strata in \cite{MZ} and more specifically in genus zero and in hyperelliptic connected components in \cite{B}.  

Usually, the study of the structure of the neighborhood of the boundary is restricted to a \emph{thick part}, where all short saddle connections are pairwise \^homo\-logous (see \cite{MS}, and also \cite{EMZ,MZ}).
Following this idea, we will consider the complement of the codimension~$1$ subset $\Delta$ of flat surfaces that admit a pair of saddle connections that have both minimal length, but which are not \^homologous.  

For a flat surface in the complement of $\Delta$, we can define the configuration of the maximal collection of \^homologous saddle connections that contains the smallest saddle connection of the surface. This defines a locally constant map outside $\Delta$ (see section \ref{CP1} for more details).

We will prove the following result.
\begin{MT} 
Let
$\mathcal{Q}_1(k_1,\ldots,k_r)$ be a stratum of quadratic differentials on $\mathbb{CP}^1$ with at most simple poles.  There is a natural bijection between the configurations of \^homologous saddle connections existing in that stratum and the connected components of $\mathcal{Q}_1(k_1,\ldots,k_r)\backslash\Delta$.
\end{MT}

We will call the connected components of $\mathcal{Q}_1(k_1,\ldots,k_r)\backslash\Delta$ the \emph{configuration domains} of the strata. 
These configuration domains might be interesting to the extend that they are ``almost'' manifolds in the following sense:

\begin{cor}
Let $\mathcal{D}$ be a configuration domain of a strata of quadratic differentials on $\mathbb{CP}^1$. If $\mathcal{D}$  admits orbifoldic points, then the corresponding configuration is symmetric and the locus of such orbifoldic points are unions of copies (or coverings) of submanifolds of smaller strata.
\end{cor}

Restricting ourselves to the neighborhood of the boundary, we show that these domains have one topological end.

\begin{prop}
Let $\mathcal{D}$ be a configuration domain of a strata of quadratic differentials on $\mathbb{CP}^1$. Let  $\mathcal{Q}_{1,\delta}(k_1,\ldots,k_r)$ be the subset of the strata corresponding to area one surfaces with at least a saddle connection of length less than $\delta$. Then $\mathcal{D}\cap \mathcal{Q}_{1,\delta}(k_1,\ldots,k_r)$ is connected for all $\delta>0$.
\end{prop}

\begin{cor}
Any stratum of quadratic differentials on $\mathbb{CP}^1$  has only one topological end.
\end{cor}

\subsubsection*{Acknowledgements} I would like to thank Anton Zorich for encouraging me to write this paper, and for many discussions.
I also thank  Erwan Lanneau and  Pascal Hubert for their useful comment.


\subsection{Basic definitions}
Here we first review standart facts about moduli spaces of quadratic differentials. We refer to \cite{Hubbard:Masur,Masur82,Veech82} for proofs and details, and to \cite{MT,Z} for general surveys.

Let  $S$  be a  compact  Riemann surface  of  genus  $g$. A  quadratic
differential $q$  on $S$ is  locally given by  $q(z)=\phi(z)dz^2$, for $(U,z)$ a local chart with  $\phi$ a meromorphic function with at most simple poles. We define the poles and zeroes of $q$ in a local chart to be the poles and zeroes of the corresponding meromorphic function $\phi$. It is easy to check that they do not depend on the choice of the local chart.  Slightly abusing notations, a marked point on the surface (\emph{resp.} a pole) will be referred to as a zero of order 0 (\emph{resp.} a zero of order $-1$).  An Abelian differential on $S$ is a holomorphic 1-form.

Outside its poles and zeros,  $q$ is locally the square of an Abelian differential.  Integrating this 1-form gives a natural atlas such that  the transition  functions are  of the kind $z\mapsto \pm z+c$. Thus $S$  inherits a flat metric with singularities, where  a zero  of order  $k\geq -1$ becomes  a conical  singularity of angle $(k+2)\pi$.  The flat metric has trivial holonomy if and only if $q$ is globally the square of any Abelian differential. If not, then the holonomy  is   $\mathbb{Z}/2\mathbb{Z}$  and  $(S,q)$  is sometimes called  a \emph{half-translation} surface since transition surfaces are either half-turns, or translations. In order to simplify the notation, we will usually denote by $S$ a surface with a flat structure.

We associate to a quadratic differential the set with multiplicity $\{k_1,\ldots,k_r\}$ of orders of its  poles and zeros.  The Gauss-Bonnet  formula asserts that  $\sum_i k_i=4g-4$. Conversely, if we fix a collection $\{k_1,\dots,k_r\}$ of integers, greater than or equal to $-1$ satisfying the previous equality, we denote by $\mathcal{Q}(k_1,\ldots,k_r)$ the (possibly empty) moduli space of quadratic differential which are not globally  squares  of Abelian differential, and which have $\{k_1,\ldots,k_r\}$ as orders of poles and zeros .  It is well known that $\mathcal{Q}(k_1,\ldots,k_r)$ is a complex analytic orbifold, which  is usually called a \emph{stratum} of the moduli space of quadratic differentials on a Riemann surface of genus $g$. We usually restrict ourselves to the subspace $Q_1(k_1,\ldots,k_r)$  of area one surfaces, where the area is given by the flat metric.  In a similar way, we denote by $\mathcal{H}_1(n_1,\dots,n_s)$ the moduli space of Abelian differentials  of area $1$ 
having zeroes of degree $\{n_1,\ldots,n_s\}$, where $n_i\geq 0$ and $\sum_{i=1}^s n_i=2g-2$.

There is a natural action of $SL_2(\mathbb{R})$ on  $\mathcal{Q}(k_1,\ldots,k_r)$ that preserve its stratification: let $(U_i,\phi_i)_{i\in I}$ is a atlas of flat coordinates of $S$, with $U_i$ open subset of $S$ and $\phi_i(U_i)\subset \mathbb{R}^2$. An atlas of $A.S$ is given by $(U_i,A\circ \phi_i)_{i\in I}$. The action of the diagonal subgroup of $SL_2(\mathbb{R})$ is called the Teichmüller geodesic flow. In order to specify notations, we denote by $g_t$ and $r_t$ the following matrix of $SL_2(\mathbb{R})$:
\begin{eqnarray*}
g_t=
\left[
\begin{array}{cc}
 e^{\frac{t}{2}} &  0    \\
 0  &    e^{-\frac{t}{2}}
\end{array}
\right]
\qquad
r_t=\left[
\begin{array}{cc}
\cos(t) & \sin(t) \\ -\sin(t) & \cos(t) 
\end{array}
\right]
\end{eqnarray*}

A saddle  connection is  a geodesic segment (or geodesic loop) joining two singularities (or a singularity to itself) with no  singularities in its interior. Even if $q$ is not  globally a square of an Abelian differential we can find a square   root   of  it   along  the   saddle connection. Integrating  it along  the  saddle  connection  we get  a complex number
(defined up  to multiplication by  $-1$). Considered as a planar vector, this complex number represents the affine holonomy vector along the saddle connection. In particular, its euclidean length is the modulus of its holonomy vector. Note that  a saddle connection persists under any small deformation of the surface. 

Local coordinates for a stratum of Abelian differential are obtained by integrating the holomorphic 1-form along a basis of the relative homology $H^1(S,{sing},\mathbb{Z})$, where $sing$ denote the set of conical singularities of $S$. Equivalently, this means that local coordinates are defined by the relative cohomology  $H_1(S,{sing},\mathbb{C})$.

Local coordinates in a stratum of quadratic differentials are obtained in the following way: one can naturally associate to a quadratic differential $(S,q)\in \mathcal{Q}(k_1,\dots,k_r)$ a double cover $p:\widehat{S}\rightarrow S$ such that $p^* q$ is the square of an Abelian differential $\omega$.  The surface $\widehat{S}$ admits a natural involution $\tau$, that induces on the relative cohomology $H_1(\widehat{S},{sing},\mathbb{C})$ an involution $\tau^*$. It decomposes $H_1(\widehat{S},{sing},\mathbb{C})$ into a invariant subspace $H_1^+(\widehat{S},{sing},\mathbb{C})$ and an anti-invariant subspace $H_1^-(\widehat{S},{sing},\mathbb{C})$. One can show that the anti-invariant subspace $H_1^-(\widehat{S},{sing},\mathbb{C})$ gives local coordinates for the stratum $\mathcal{Q}(k_1,\ldots,k_r)$.

\subsection{\^Homologous saddle connections}\label{homolo}

Let $S\in \mathcal{Q}(k_1,\dots,k_r)$ be a flat surface and denote by  $p:\widehat{S}\rightarrow S$ its canonical double cover and $\tau$ its corresponding involution. Let $\Sigma$ be the set of singularities of $S$ and $\widehat{\Sigma}=p^{-1}(\Sigma)$.

To an oriented saddle connection $\gamma$ on $S$, we can associate $\gamma_1$ and $\gamma_2$ its preimages by $p$. If the relative cycle $[\gamma_1]$ satisfies $[\gamma_1]=-[\gamma_2] \in H_1(\widehat{S},\widehat{\Sigma},\mathbb{Z})$, then we define $[\hat{\gamma}]=[\gamma_1]$. Otherwise, we define  $[\hat{\gamma}]=[\gamma_1]-[\gamma_2]$. Note that in all cases, the cycle $[\hat{\gamma}]$ is anti-invariant with respect to the involution~$\tau$.

\begin{defs} 
Two saddle connections $\gamma$ and $\gamma^{\prime}$ are \^homologous if $[\hat{\gamma}]=\pm [\hat{\gamma^{\prime}}]$.
\end{defs}

\begin{figure}[htb]
\begin{center}
\input{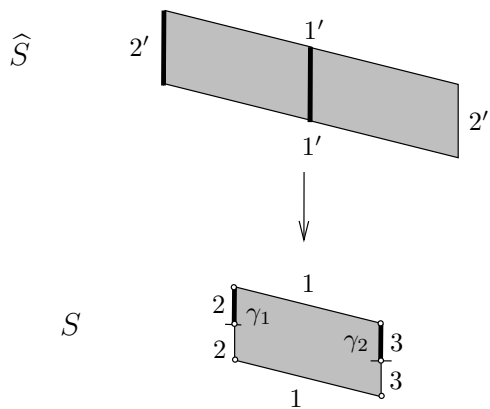}
\caption{An unfolded flat surface $S$ with two \^homo\-logous saddle connections $\gamma_1$ and $\gamma_2$.}
\label{exemple_trivial}
\end{center}
\end{figure}

\begin{Example}
Consider the flat surface $S\in \mathcal{Q}(-1,-1,-1,-1)$ given in Figure \ref{exemple_trivial} (a ``pillowcase''), it is easy to check from the definition that $\gamma_1$ and $\gamma_2$ are \^homo\-logous since the corresponding cycles for the double cover $\widehat{S}$ are homologous.
\end{Example}

\begin{figure}[htp]
\begin{center}
\input{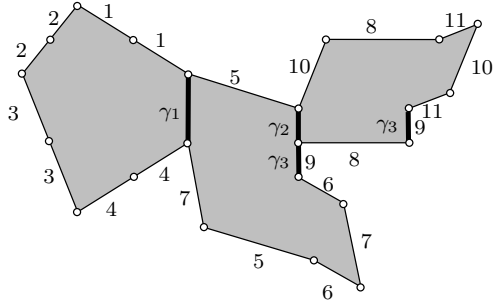}
\caption{Unfolded flat surface with three \^homo\-logous saddle connections $\gamma_1$, $\gamma_2$, and $\gamma_3$.}
\label{exempleplat}
\end{center}
\end{figure}

\begin{Example}
Consider the flat surface given in Figure \ref{exempleplat}, the reader can check that the saddle connections $\gamma_1$, $\gamma_2$ and $\gamma_3$ are pairwise \^homo\-logous.
\end{Example}

The following theorem is due to Masur and Zorich \cite{MZ}. It gives in particular a simple geometric criteria for deciding whether two saddle connections are \^homo\-logous.  We give in the appendix an alternative proof.

\begin{NNths}[H. Masur; A. Zorich] 
Consider two distinct saddle connections $\gamma,\gamma^{\prime}$ on a half-translation surface. The following assertions are equivalent:
 \begin{itemize} 
 \item  The two saddle connections $\gamma$ and $\gamma^{\prime}$ are \^homologous.  
\item The ratio of their length is constant under any small deformation of the surface inside the ambient stratum.  
\item
They  have no interior intersection and one of the connected  component  of
$S\backslash\{\gamma\cup \gamma^{\prime}\}$ has  trivial  linear holonomy.
\end{itemize} 
Furthermore, if $\gamma$ and $\gamma^{\prime}$ are \^homologous, then the ratio of their length belongs to $\{1/2,1,2\}$ and they are parallel.
\end{NNths}

A saddle connection $\gamma_1$ will be called \emph{simple} if they are no other saddle connections \^homologous to $\gamma_1$.  Now we consider a set of \^homologous saddle connections $\gamma=\{\gamma_1, \dots,\gamma_s\}$ on a flat surface $S$. 
Slightly abusing notation, we will denote by $S\backslash   \gamma$ the subset $S\backslash\bigl( \cup_{i=1}^s \gamma_i \bigr)$. This subset is a finite union of connected
half-translation surfaces with boundary. We define a graph $\Gamma(S,\gamma)$  called  the  graph  of connected  components in  the following  way  (see  \cite{MZ}):  the vertices  are  the  connected components   of  $S\backslash   \gamma$,  labelled   as  ``$\circ$''   if  the corresponding surface is a cylinder,  as ``$+$'' if it has trivial holonomy (but is  not a cylinder), and  as ``$-$'' if it  has non-trivial holonomy.
The  edges are  given  by  the saddle  connections  in $\gamma$.  Each
$\gamma_i$ is  on the  boundary of one  or two connected  components of
$S\backslash  \gamma$. In  the  first case  it  becomes an  edge joining  the
corresponding vertex to itself. In the second case, it becomes an edge
joining the two corresponding vertices. 

Each connected  components of $S\backslash \gamma$ is a non-compact surface but can be naturally compactified (for example considering the distance induced by the flat metric on a connected component of $S\backslash \gamma$, and the corresponding completion).  We denote this compactification by $S_j$. We warn the reader that $S_j$ might  differ from the  closure of the component in the surface $S$: for example, if $\gamma_i$ is on the boundary of just one connected  component $S_j$ of $S\backslash\gamma$ , then the compactification of $S_j$ contains  two copies of $\gamma_i$ in its boundary,  while  in  the  closure  of  $S_j$  these  two copies  are identified.  The  boundary of each  $S_i$ is a union of saddle connections; it has one or  several connected components.  Each of them is homeomorphic to $\mathbb{S}^1$ and therefore the orientation of $S$ defines a cyclic order in the set of boundary saddle connections.  Each consecutive pair of saddle connections for that cyclic order defines a \emph{boundary singularity} with an associated angle which is a integer multiple of $\pi$ (because the boundary saddle connections are parallel). The surface with boundary $S_i$ might have singularities in its interior. We call them \emph{interior singularities}.

\begin{defs} 
Let $\gamma=\{\gamma_1,\dots,\gamma_r\}$ be a maximal collection of \^ho\-mo\-lo\-gous saddle connections on a flat surface. A \emph{configuration} is the following combinatorial data: 
\begin{itemize} 
\item The graph $\Gamma(S,\gamma)$ 
\item For each vertex of this graph, a permutation of the edges adjacent to the vertex (encoding the cyclic order of the saddle connections on each connected component of the boundary of the $S_i$).  
\item For each pair of consecutive elements in that cyclic order, the angle between the two corresponding saddle connections. 
\item For each $S_i$, a collection of integers that are the orders of the interior singularities of $S_i$.  
\end{itemize} 
\end{defs}

We refer to \cite{MZ} for a more detailed definition of a configuration (see also \cite{B}).

\subsection{Neighborhood of the boundary, thick-thin decomposition}

 For any compact subset $K$ of a stratum, there exists a constant $c_{K}$ such that the length of any saddle connection of any surface in $K$ is greater than $c_{K}$.  Therefore, we can define the $\delta$-neighborhood of the boundary of the strata to be the subset of area $1$ surfaces that admit a saddle connection of length less than $\delta$.

According to Masur and Smillie \cite{MS}, one can decompose the  $\delta$-neigborhood of the boundary of a stratum into a  \emph{thin part} (of negligibly small measure) and a \emph{thick part}. The thin part being for example  the subset of surfaces with a pair of non\^homo\-logous saddle connections of length respectively less than $\delta$ and $N\delta$, for some fixed $N\geq 1$ (the decomposition depends on the choice of $N$).  We also refer to \cite{EMZ} for the case of Abelian differentials and  to \cite{MZ} for the case of quadratic differentials.

Let $N\geq 1$, we consider  $Q^{N}(k_1,k_2,\dots,k_r)$  the subset of flat surfaces  such that, if $\gamma_1$ is the shortest saddle connection and $\gamma_1^{\prime}$ is another saddle connection non\^homologous to $\gamma_1$, then $|\gamma_1^{\prime}|> N|\gamma_1|$. Similarly, we define $Q_1^{N}(k_1,k_2,\dots,k_r)$ the intersection of $Q^{N}(k_1,k_2,\dots,k_r)$ with the subset of area $1$ flat surfaces.

 For any surface in $Q^{N}(k_1,k_2,\dots,k_r)$, we can define a maximal collection $\mathcal{F}$ of \^homologous saddle connections that contains the smallest one. This is well defined because if there exists two smallest saddle connections, they are necessary \^homologous. We will show in section \ref{CP1} the associated configuration defines a locally constant map from $Q_1^{N}(k_1,k_2,\dots,k_r)$ to the space of configurations. This leads to the following definition:
 
\begin{defs} 
A \emph{configuration domain} of $\mathcal{Q}_1(k_1,\dots,k_r)$ is a connected component of $Q_1^{N}(k_1,\ldots,k_r)$.  
\end{defs}

\begin{rem}
The previous definition of a configuration domain is a little more general than the one stated in the introduction that corresponds to the case $N=1$.

\end{rem}

\begin{defs}
An \emph{end} of a space $W$ is a function  \[ \epsilon: \{K,\ K\subset W \textrm{ is compact}\} \rightarrow \{X,\ X\subset W\}\] such that:
\begin{itemize}
\item $\epsilon(K)$ is a (unbounded) component of $W\backslash K$ for each $K$
\item if $K\subset L$, then $\epsilon(L)\subset \epsilon(K)$.
\end{itemize}
\end{defs}

\begin{NNprop}
If $W$ is $\sigma-$compact, then the number of ends of $W$ is the maximal number of unbounded components of $W\backslash K$, for $K$ compact.
\end{NNprop}

We refer to \cite{HR} for more details on the ends of a space.

\subsection{Example on the moduli space of flat torus}

If $T$ is a flat torus (\emph{i.e.} a Riemann surface with an Abelian differential $\omega$), then, up to rescaling $\omega$, we can assume that  the holonomy vector of the shortest geodesic is $1$. Then, choosing a second smallest non horizontal geodesic with a good choice of its orientation, this defines a complex number $z=x+iy$, with  $y>0$ , $-1/2\leq x\leq 1/2$ and $|z|\geq 1$. The corresponding domain $\mathcal{D}$ in $\mathbb{C}$ is a fundamental domain of $\mathbb{H}/{SL_2(\mathbb{Z})}$.

 It is well know that this defines an map from  the moduli space of flat torus with trivial holonomy (\emph{i.e.} $\mathcal{H}(\emptyset)$), to $\mathbb{H}/{SL_2(\mathbb{Z})}$ which is a  bundle, with $\mathbb{C}^*$ as fiber. Orbifoldic points of $\mathcal{H}(\emptyset)$ are over the complex number $z_1=i$ and $z_2=\frac{1+i\sqrt{3}}{2}$. They correspond to Abelian differential on torus obtained by identifying the opposite sides of a square, or a regular hexagon. 

\begin{figure}[htbp]
   \begin{center}
      \input{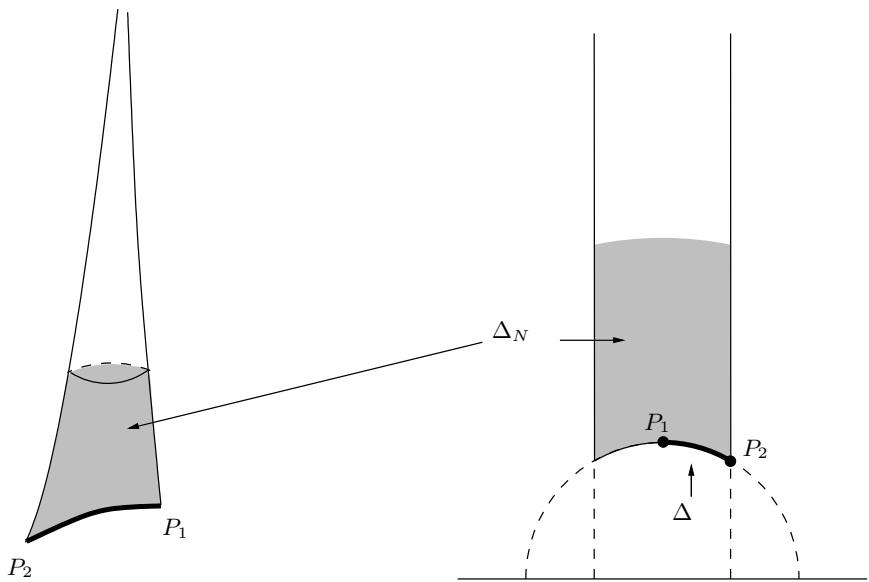}
      \caption{Configuration domain in $\mathcal{H}(\emptyset)$.}
      \label{config:tore}
   \end{center}
\end{figure}

Now with this representation, $\mathcal{H}^N(\emptyset)$ is obtained by restricting ourselves to the subdomain $\mathcal{D}^N=\mathcal{D}\cup \{z, |z|> N\}$ (see Figure \ref{config:tore}). This subdomain contains neither $z_1$ nor $z_2$, so $\mathcal{H}^N(\emptyset)$ is a manifold. In the extreme case $N=1$, the codimension one subset $\Delta$ is 
an arc joining $z_1$ to $z_2$.

\subsection{Reader's guide} 
Now we sketch the proof of the Main Theorem.
\begin{enumerate} 
 \item We first prove the theorem for the case of configuration domains defined by a simple saddle connection (we will refer to these configuration domains as \emph{simple}). We will explain how we can shrink a simple saddle connection, when its length is small enough (therefore, describe the structure of the stratum in a neighborhood of an adjacent one).  This is done in section~\ref{simple_cusps}.

There is one easy case, when the shrinking process is done by local and canonical surgeries. The other case involves some non-local surgeries (hole transport) that depend on a choice of a path. We will have to  describe the dependence of the choice of the path. More details on these surgeries appears in section~\ref{hole_transport}.

\item The list of configurations was established by the author in  \cite{B}. The second step of the proof is to consider each configuration and to show that the subset of surface associated to this configuration is connected. This will be done in section \ref{CP1} and will use the ``simple case''.
\end{enumerate}


\section{Families of quadratic differentials defined by an involution} \label{permut}

Consider a polygon whose sides come by pairs, and such that, for each pair, the corresponding sides are parallel and have the same length. Then identifying these pair of sides by appropriate isometries, this gives a flat surface. In this section we show that any flat surface can arise from such a polygon and give an explicit construction. We end by a technical lemma that will be one of the key arguments of Theorem \ref{cusp_elementaire}.

\begin{figure}[htbp]
   \begin{center}
      \input{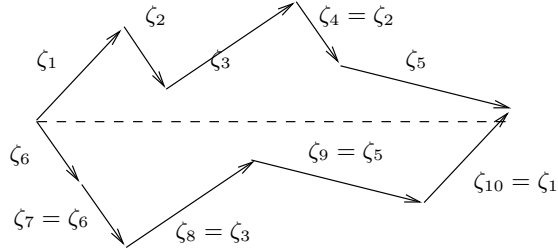}
      \caption{Flat surface unfolded into a polygon.}
      \label{polygone}
   \end{center}
\end{figure}

\subsection{Constructions of a flat surface}

Let $\sigma$ be an involution of $\left\{ 1,\dots,l+m \right\}$, without fixed points.

We denote by $Q_{\sigma,l}$ the set  $\zeta=(\zeta_1,\dots,\zeta_{l+m}) \in \mathbb{C}^{l+m}$ such that:
\begin{enumerate} 
\item $\forall i \quad \zeta_i=\zeta_{\sigma(i)}$
 %
%
\item $\forall i \quad Re(\zeta_i)>0$.  
\item $\forall 1\leq i \leq l-1 \quad Im(\sum_{k\leq i} \zeta_k)>0$ 
\item $\forall 1\leq j \leq m-1 \quad Im(\sum_{1\leq k\leq j} \zeta_{l+k})<0$
\item $\sum_{k\leq l} \zeta_k = \sum_{1\leq k\leq m}\zeta_{l+k}$.
\end{enumerate} 

Note that $Q_{\sigma,l}$ is convex and might be empty for some $\sigma$. 
Now we will construct a map $ZR$ from $Q_{\sigma,l}$ to the moduli space of quadratic differentials. Slightly abusing conventional terminology, we will call a  surface in $ZR(Q_{\sigma,l})$ a \emph{suspension} over $(\sigma,l)$, and a vector in   $Q_{\sigma,l}$ is then a \emph{suspension data}.

Furthermore, since $Q_{\sigma,l}$ is convex, the connected component of the stratum is uniquely determined by ($\sigma,l)$.

\subsubsection*{Easy case}
Now we consider a broken line $L_1$ whose edge number $i$ ($1\leq i \leq l$) is represented  by the complex number $\zeta_i$. Then we consider a second broken line $L_2$ which starts from the same point, and whose edge number $j$ ($1 \leq j\leq m$) is represented by $\zeta_{l+j}$. The last condition implies that these two lines also end at the same point. 
If that they have \emph{no other intersection points}, then they form a polygon (see Figure \ref{polygone}). The sides of the polygon, enumerated by indices of the corresponding complex number, naturally come by pairs according to the involution $\sigma$. Gluing these pair of sides by isometries respecting the natural orientation of the polygon, this construction defines a flat surface which have trivial or non-trivial holonomy.


\subsubsection*{First return map on a horizontal segment}
Let $S$ be a flat surface and $X$ be a horizontal segment with a choice of a positive vertical direction (or equivalently, a choice of left and right ends). We consider the first return map $T_1:X\rightarrow X$ of the vertical geodesic flow in the positive direction. Any infinite vertical geodesic starting from $X$ will intersect $X$ again. Therefore, the map $T_1$ is well defined outside a finite number of points that correspond to vertical geodesics that stop at a singularity before intersecting again the interval $X$. This set $X\backslash \{sing\}$ is a finite union  $X_1,\ldots,X_{l}$ of open intervals and the restriction of $T_1$ on each $X_i$  is of the kind $x\mapsto \pm x+c_i$. The time of first return of the geodesic flow is constant along each $X_i$. Similarly, we define $T_2$ to be the first return map of the  vertical flow in the negative direction and denote by $X_{l+1},\ldots,X_{l +m}$ the corresponding intervals. Remark that for $i\leq l$ (resp. $i>l$) , $T_1 (X_i)=X_j$  (resp.
$T_2(X_i)=X_j$) for some $1\leq j\leq l+m$. Therefore, $ (T_1,T_2)$ induce a permutation $\sigma_X$ of $\{1,l+m\}$, and it is easy to check that $\sigma_X$ is an involution without fixed points. When $S$ is a translation surface, $T_2=T_1^{-1}$ and $T_1$ is called an \emph{interval exchange transformation}.

If $S\in ZR(Q_{\sigma,l})$, constructed as previously, we choose $X$ to be the horizontal line whose left end is the starting point of the broken lines, and of length $Re(\sum_{k\leq l} \zeta_k)$ . Then it is easy to check that $\sigma_X=\sigma$. 

\subsubsection*{Veech zippered rectange construction}
\begin{figure}[htbp]
   \begin{center}
     \input{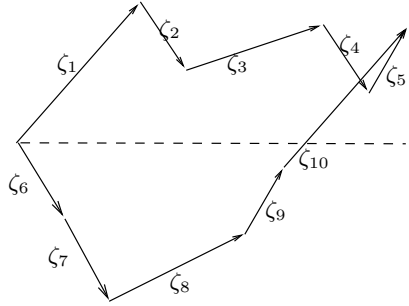}
     \caption{Suspension data that does not give a ``suitable'' polygon.}
     \label{figure:wrong:polygon}
   \end{center}
\end{figure}

The broken lines $L_1$ and $L_2$ might intersect at other points (see Figure \ref{figure:wrong:polygon}).
However, we can still define a flat surface by using an analogous construction as the well known zippered rectangles construction due to Veech. We give a description of this construction and refer to \cite{Veech82,Y}
for the case of Abelian differentials. This construction is very similar to the usual one, although its precise description is quite technical. Still, for completeness, we give an equivalent but rather implicit formulation.

\begin{figure}[htbp]
   \begin{center}
     \input{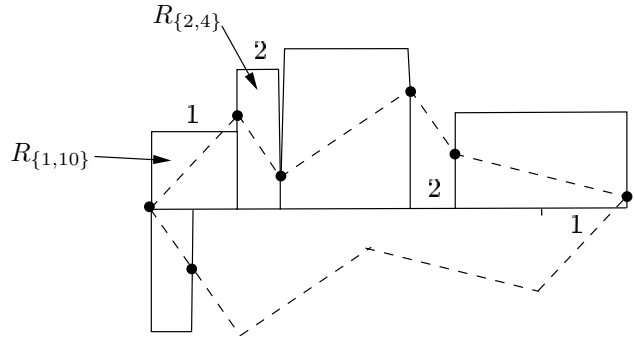}
     \caption{Zippered rectangle construction, for the case the flat surface of Figure \ref{polygone}.}
     \label{polygtozip}
   \end{center}
\end{figure}

We first consider the previous case when $L_1$ and $L_2$ define an acceptable polygon. For each pair of interval $X_i, X_{\sigma(i)}$ on $X$, the return time $h_i=h_{\sigma(i)}$ of the vertical flow starting from $x\in X_i$ and returning in $y\in X_{\sigma(i)}$ is constant. This value depends only on $(\sigma,l)$ and on the imaginary part of  $\zeta$. For each pair $\alpha=\{i,\sigma(i)\}$ there is a natural embedding of the open rectangle $R_{\alpha}=(0,Re(\zeta_i))\times(0,h_i)$ into the flat surface $S$ (see Figure \ref{polygtozip}). For each $R_\alpha$, we glue a horizontal side to $X_i$ and the other to $X_{\sigma(i)}$. The surface $S$ is then obtained after suitable identifications of the vertical sides of the the rectangles $\{R_\alpha\}$. These vertical identifications  only depend on $(\sigma,l)$ and on the imaginary part of $\zeta$. 

For the general case, we construct  the rectangles $R_{\alpha}$ by using the same formulas. Identifications for the horizontal sides are staightforward. Identifications for the vertical sides  do not depends on the horizontal parameters, and will be the same as for a suspension data $\zeta^{\prime}$ that have the same imaginary part as $\zeta$, but that correspond to a suitable polygon. This will be well defined after the following lemma.

\begin{lem}
Let $\zeta$ be a collection of complex numbers in  $Q_{\sigma,l}$ then there exists $\zeta^{\prime} \in Q_{\sigma,l}$ with the same imaginary part as $\zeta$, that defines a suitable polygon.
\end{lem}
\begin{proof}
We can assume that $\sum_{k=1}^{l} Im(\zeta_k)>0$ (the negative case is analogous and there is nothing to prove when the sum is zero).
It is clear that $\sigma(l+m)\neq l$ otherwise there would be no possible suspension data. If $\sigma(l+m)<l$, then we can shorten the real part of  $\zeta_{l+m}$ and of $\zeta_{\sigma(l+m)}$, keeping conditions (1)---(5) satisfied, and get a suspension data $\zeta^{\prime}$ with the same imaginary part as $\zeta$, and such that $Re(\zeta^{\prime}_{l+m})<Re(\zeta^{\prime}_l)$. This last condition implies that $\zeta^{\prime}$ defines a suitable polygon.

If $\sigma(l+m)>l$, then condition $(5)$ might imply that $Re(\zeta_{l+m})$ is necessary bigger than $Re(\zeta_l)$. However, we can still change $\zeta$ into a suspension data $\zeta^{\prime}$, with same imaginary part, and such that $Re(\zeta^{\prime}_{l+m})$ is very close to $Re(\zeta^{\prime}_l)$. In that case, $\zeta^{\prime}$ also defines a suitable polygon.
\end{proof}

\subsection{The converse: construction of suspension data from a flat surface}


Now we give a sufficient condition for a surface to be in some $Q_{\sigma,l}$. 
Note that an analogous construction for hyperelliptic flat surfaces has been done in \cite{Veech:hyp}.

\begin{figure}[htbp]
   \begin{center}
      \input{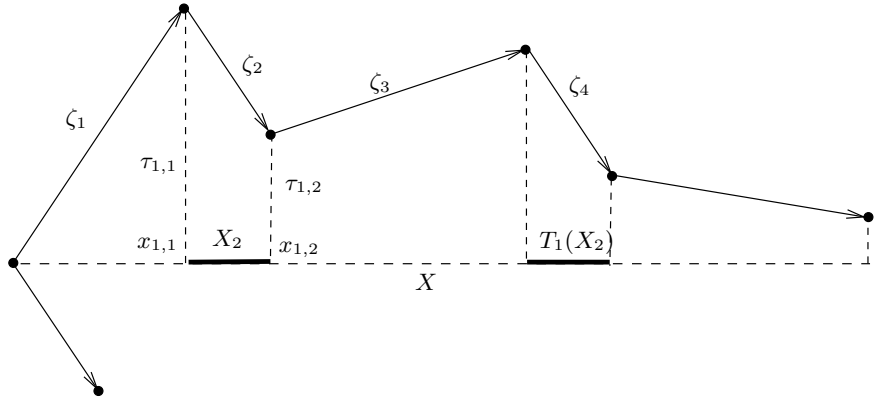}
      \caption{Construction of a polygon from a surface.}
   \end{center}
\end{figure}


\begin{prop} \label{suspension}
Let $S$ be a flat surface with no vertical saddle connection. There exists an involution $\sigma$ and an integer $l$ such that $S\in ZR(Q_{\sigma,l})$.
\end{prop}

\begin{figure}[htbp]
   \begin{center}
      \input{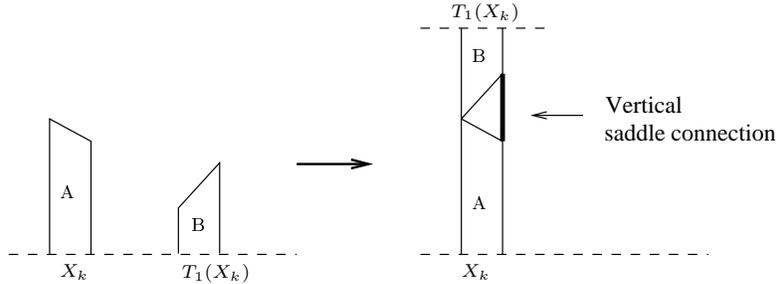}
      \caption{The complex numbers $\zeta_k$ and $\zeta_{\sigma_X(k)}$ are necessary equal.}
         \label{sc_verticale}
   \end{center}
\end{figure}

\begin{proof} 
Let $X$ be a horizontal segment whose left end is a singularity. Up to cutting $X$ on the right, we can assume that the vertical geodesic starting from its right end hits a singularity before meeting $X$ again.
	
Let $x_{1,1}<\ldots<x_{1,l-1}$ be the points of discontinuity of $T_1$ and $(x_{1,0},x_{1,l})$ be the end points of $X$. For each $k$, there exists $\tau_{1,k}>0$ such that the vertical geodesic starting from $x_{1,k}$ in the positive direction stops at a singularity at time $\tau_{1,k}$ (here $\tau_{1,0}=0$, since by convention $x_{1,0}$ is located at a singularity). Then for $k \geq 1$ we define $\zeta_k: (x_{1,k}-x_{1,k-1})+i(\tau_{1,k}-\tau_{1,k-1})$.  Now we perform a similar construction for the vertical flow in the negative direction: let $x_{2,1}<\ldots<x_{2,m-1}$ be the points of discontinuity of $T_2$ and $(x_{2,0},x_{2,m})$ be the extremities of $X$. For each $k\notin\{0,m\}$, the vertical geodesic starting from $x_{2,k}$ in the positive direction stops at a singularity at time $\tau_{2,k}<0$ (here again $\tau_{2,0}=0$ and $\tau_{2,l}>0$).  For $1\leq k\leq m$, we define $\zeta_{k+l}: (x_{2,k}-x_{2,k-1})+i (\tau_{2,k}-\tau_ {2,k-1})$. So, we have a collection of complex numbers $\zeta_{l+1},\ldots,\zeta_{m +l}$ that defines a polygon $\mathcal{P}$.

We have always $Re(\zeta_k)=Re(\zeta_{\sigma_X(k)})=|X_k|$.  Let $1\leq k\leq l$. If $\sigma_X(k)\leq l$, then  $\tau_{1,k-1}+\tau_{1,\sigma_X(k)}=\tau_{1,k}+\tau_{1,\sigma_X(k)-1}=h_k$ (with $h_k$ the time of first return to $X$ of the vertical geodesic flow starting from the subinterval $X_k$), otherwise there would exist a vertical saddle connection (see figure \ref{sc_verticale}). So  $Im(\zeta_k)=Im(\zeta_{\sigma_X(k)})$.  The other cases are analogous. 
Thus $\zeta$ is a suspension data, and $ZR(\zeta)$ is isometric to $S$.
\end{proof}

\begin{rem}
In the previous construction, the suspension data constructed does not necessary give a ``suitable'' polygon. However, we can show that by choosing carefully a subinterval $X^{\prime}$ of $X$, the construction will give a true polygon. Since for any surface, we can find a direction with no saddle connection, we can conclude that any surface can be unfolded into a polygon as in Figure \ref{polygone}, up to rotating that polygon.
\end{rem}


\subsection{A technical lemma}

The following lemma is a technical lemma that will be needed in section $\ref{non_loc}$. It can be skipped in a first reading.  We previously showed that a surface with no vertical saddle connection belongs to some $ZR(Q_{\sigma,l})$. Furthermore, the corresponding pair $(\sigma,l)$ is completely defined by a first return map of the vertical flow on a well chosen horizontal segment.

We define the set $Q_{\sigma,l}^{\prime}$ defined in a similar way as
$Q_{\sigma,l}$, but here we replace condition $2$  by the following two conditions:
 \begin{enumerate} 
 \item[$(2)$] $\forall i\notin \{1,\sigma(1)\} \quad Re(\zeta_i)>0$.  
 \item[$(2^{\prime})$] $Re(\zeta_1)=0$.  
 \end{enumerate}
In other words, the first vector of the top broken line is now vertical and no other vector is vertical except the other one of the corresponding pair.
Then we define in a very similar way a map $ZR^{\prime}$ from $Q_{\sigma,l}^{\prime}$ to a stratum of the moduli space of quadratic differentials.

Note that the subset $Q_{\sigma,l}^{\prime}$ is convex.

\begin{lem} \label{key_fam} 
 Let $S$ be a flat surface with a unique vertical saddle connection joining two singularities $P_1$ and $P_2$. Let $X$ be a horizontal segment whose left end is $P_1$, and such that the vertical geodesic starting from its left end is the unique vertical saddle connection joining $P_1$ to $P_2$. There exists $(\sigma,l)$ that depends only on the first return maps on $X$ of the vertical flow and on the degree of $P_2$, such that $S\in ZR^{\prime}(Q_{\sigma,l}^{\prime})$.  
\end{lem}

\begin{proof} 
 We define as in Proposition \ref{suspension} the $x_{i,j} , \tau_{i,j}$ and $\zeta_j$, with the slight difference that now, $\tau_{1,0}>0$. Now, because there exists only one vertical saddle connection, the same argument as before says that there exists at most one unordered pair $\{\zeta_{i_0}, \zeta_{\sigma(i_0)}\}$ such that $\zeta_{i_0}\neq \zeta_{\sigma(i_0)}$.
If this pair doesn't exists, then the union of the vertical geodesics starting from $X$ would be a strict subset of $S$, with boundary the unique vertical saddle connection. Therefore, we would have $P_1=P_2$, contradicting the hypothesis.

Now we glue on the polygon $\mathcal{P}$ an Euclidean triangle of sides given by $\{\zeta_i, \zeta_{\sigma(i)},i\tau_{1,0}\}$, and we get a new polygon. The sides of this polygon appear in pairs that are parallel and of the same length. We can therefore glue this pair and get a flat surface. By construction, we get a surface isometric to $S$, and so $S$ belongs to some $ZR^{\prime}(Q_{\tilde{\sigma},l}^{\prime})$. The permutation $\tilde {\sigma}$ is easily constructed from $\sigma$ as soon as we know $i_0$.  This value is obtained by the following way: we start from the vertical saddle connection, close to the singularity $P_2$. Then, we turn around $P_2$ counterclockwise. Each half-turn is easily described in terms of the permutation $\sigma$. Then after performing $k_2+2$ half-turns, we must arrive again on the vertical saddle connection. This gives us the value of $i_0$.  
%

\end{proof}

\section{Hole transport} \label{hole_transport}

Hole transport is a  surgery used in \cite{MZ} to show the existence of some configurations and especially to break an even singularity to a pair of odd ones. It was defined along a simple path transverse to the vertical foliation. In this part, we generalize this construction to a larger class of paths and show that breaking a zero using that procedure  does not depend on small perturbations of the path.

Hole transport also appears in \cite{EMZ} in the computation of the Siegel-Veech constants for the moduli space of Abelian differentials. 
This improved surgery, and ``dependence properties'' that are  Corollary \ref{petites:variations} and Lemma \ref{independant} are a necessary toolkit for the computayion of these these Siegel-Veech constants for the case of quadratic differentials.

\begin{defs} A hole is a connected component of the boundary of a flat surface given by a single saddle connection. The saddle connection bounds a singularity. If this singularity has angle $3\pi$, this hole is said to be simple.  \end{defs} 

\begin{conv} 
We will always assume that the saddle connection defining the hole is vertical \end{conv}

A simple hole $\tau$ has a natural orientation given by the orientation of the underlying Riemann surface. In a neighborhood of the hole, the flat metric has trivial holonomy and therefore $q$ is locally the square of an Abelian differential.

\begin{conv} \label{orient_trou} 
When defining the surgeries around a simple hole using flat coordinates, we will assume (unless explicit warning) that the flat coordinates come from a local square root $\omega$ of $q$, such that that $\int_{\tau}dz\in i \mathbb{R}^+$.  \end{conv}

\begin{rem} 
Under convention \ref{orient_trou}, we may speak of the \emph{left} or the \emph {right} direction in a neighborhood of a simple hole. Note that there exists two horizontal geodesics starting from the singularity of  and going to the right, and only one starting from the singularity and going to the left.
\end{rem}

\begin{figure}[htpb]
\begin{center}
\input{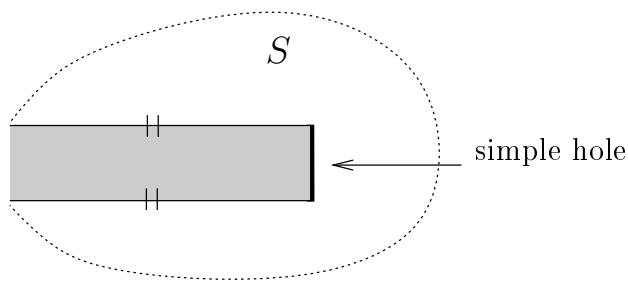}
\caption{A hole in flat coordinates.}
\end{center}
\end{figure}


\subsection{Parallelogram constructions}
We first describe the three basic surgeries on the surface that allow us to transport a simple hole along a segment. Consider a simple hole $\tau$ and chose flat coordinates in a neighborhood of the hole that satisfy convention 2.  We consider a vector $v$ such that $Re(dz(v))>0$ (\emph{i.e.} the vector $v$ goes ``to the right'' in our flat coordinates).  Consider the domain $\Omega$ obtained as the union of geodesics of length $|v|$, starting at a point of $\tau$ with direction $v$. When $\Omega$ is an embedded parallelogram, we can remove it and glue together by translation the two sides parallel to $v$.  Here we have transported the simple hole by the vector $v$. Note that the area changes under this construction.

When $Re(dz(v))<0$, this construction (removing a parallelogram) cannot work. The singularity is the unique point of the boundary that can be the starting point of a geodesic of direction $v$.  Now from the corresponding geodesic, we perform the reverse construction with respect to the previous one: we cut the surface along a segment of length $v$  and paste in a parallelogram. By means of this construction we transport the hole along the vector $v$. 

When $Re(dz(v))=0$, we consider a geodesic segment of direction $v$ starting from the singularity, and cut the surface along the segment, then glue it with a shift (``Earthquake construction'').

\begin{figure}[htpb]
\begin{center}
\input{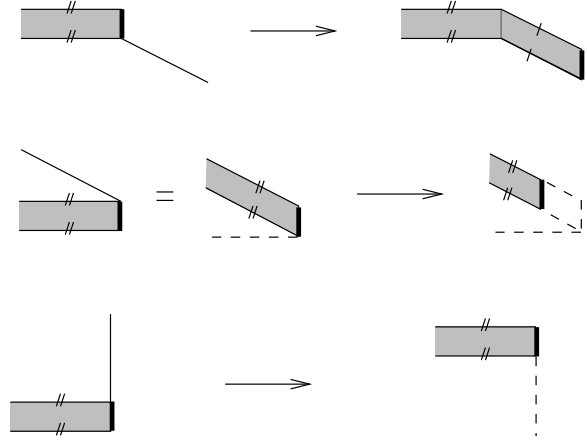}
\caption{Parallelogram constructions.}
\end{center}
\end{figure}

There is an easy way to create a pair of holes in a compact flat surface: we consider a geodesic segment imbedded in the surface, we cut the surface along that segment and paste in a parallelogram as in the previous construction.
We get parallel holes of the same length (but with opposite orientation).  Note that we can assume that the length of these holes is arbitrary small. In a similar way, we can create a pair of holes by removing a parallelogram.

\subsection{Transport along a piecewise geodesic path}
Now we consider a piecewise geodesic simple path $\gamma=\gamma_1\ldots\gamma_n$ with edges represented by the vectors $v_1,v_2,\dots,v_n$. We assume for simplicity that neither of $v_i$ is vertical.  The spirit is to transport the hole by iterating the previous constructions. We make the hole to ``follow the path'' $\gamma$ in the following way (under convention \ref{orient_trou}):

\begin{itemize} 
 \item At step number $i$, we ask that the geodesic $\gamma_i$ starts from the singularity of the hole.  \item When $Re(dz(v_i))>0$, we ask $\gamma_i$ to be the bottom of the parallelogram $\Omega$ defined in the previous construction.  \end{itemize}

\begin{figure}[htbp]
\begin{center}
\input{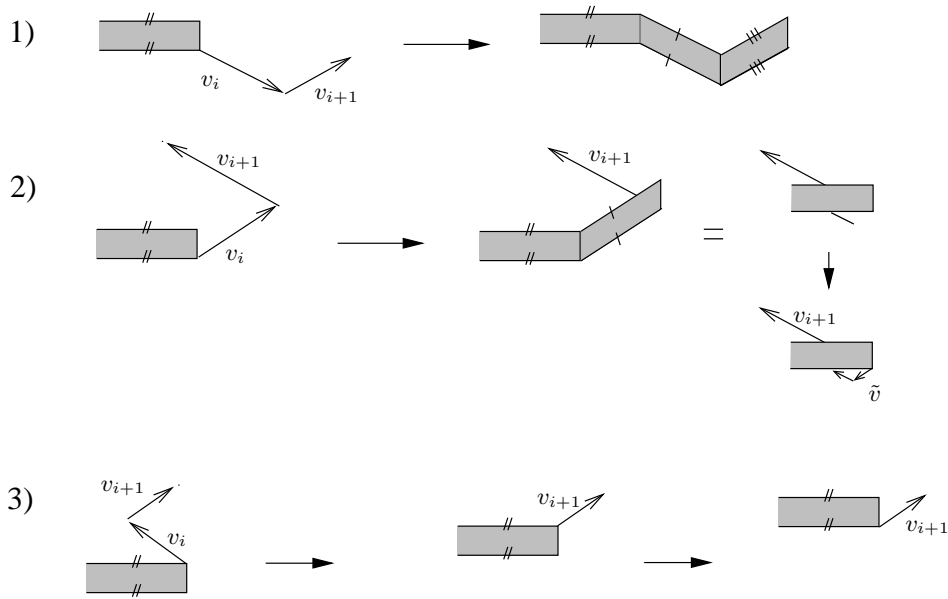}
\caption{Hole transport along a piecewise geodesic curve.}
\label{iterating}
\end{center}
\end{figure}

Naive iteration does not necessary preserve these conditions.  The surgery can indeed disconnect the path  but then we can always reconnect $\gamma$ by adding a geodesic segment. If the first condition is satisfied, but not the second, we can add a surgery along a vertical segment of the size of the hole to fulfill it.  We just have to check that each iteration between two consecutive segments of the initial path can be done in a finite number of steps, see Figure \ref{iterating}.

\begin{enumerate}
\item
If  $Re(dz(v_i))$ and $Re(dz(v_{i+1}))$ have the same sign, then as soon as both transports are successively possible, our two conditions keep being fulfilled.

\item
If $Re(dz(v_i))>0$ and $Re(dz(v_{i+1}))<0$, and if $(v_i,v_{i+1})$ is positively oriented, the surgery with $v_i$ disconnect the path, and we must add a new segment $\tilde{v}$, but then $Re(\tilde{v})$ and $Re(v_{i+1})$ are both negative, therefore, we can iterate the surgery keeping the two conditions fulfilled.

\item
If $Re(dz(v_i))<0$ and $Re(dz(v_{i+1}))>0$, and if $(v_i,v_{i+1})$ is negatively oriented, we must add a surgery along a vertical segment to fulfill the second condition.

\item
It is an easy exercise to check that for any other configuration of $(v_i,v_{i+1})$, the direct iteration of the elementary surgeries works.
\end{enumerate}

Of course, in the process we have just described, we implicitly assumed that at each step, the condition imposed for the basic surgeries (\emph{i.e.} the parallelogram must be  imbedded in the surface) is fulfilled. But considering any compact piecewise geodesic path, the process will be well defined as soon as the hole is small enough.

\begin{rem}
We can also define hole transport along a piecewise geodesic path that have self intersections. Here hole transport will disconnect the path at each intersections, but we can easily reconnect it and hole transport  also ends in a finite number of steps. We will not need  hole transport along such paths.
\end{rem}

\subsection{Application: breaking up an even singularity} \label{break_non_loc}
We consider a singularity $P$ of order $k=k_1+k_2$. When $k_1$ and $k_2$ are not both odd, there is a local surgery that  continuously break this singularity into pair of singularities of order $k_1$ and $k_2$ (see section \ref{break_loc}). When $k_1$ and $k_2$ are both odd, this local surgery fails. Following \cite{MZ} we use hole transport instead.

\begin{figure}[htpb]
\begin{center}
\input{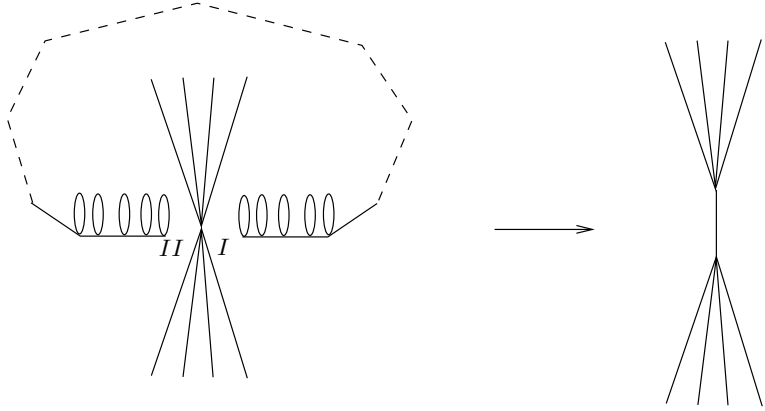}
\caption{Breaking a singularity.}
\end{center}
\end{figure}

Consider a pair $(I,II)$ of sectors of angle $\pi$ in a small neighborhood of $P$, and such that the image of the first one by a rotation of $(k_2+1)\pi$ is the second sector. Now let $\gamma$ be a simple broken line that starts and ends at $P$, and such that its first segment belongs to sector $I$ and its last segment belongs to sector $II$. We ask that parallel transport along $\gamma$ to be $\mathbb{Z}/2\mathbb{Z}$ (this has sense because $k$ is even, so $P$ admits a parallel vector field in its neighborhood).

Then, we create a pair of holes by cutting the first segment and pasting in a parallelogram. Denote by $\varepsilon$ the length of these holes. One hole is attached to the singularity. The other one is a simple hole.  We can transport it along $\gamma$, to the sector $II$. Then gluing the holes together, we get a singular surface with a pair of conical singularities that are glued together. If we desingularise the surface, we get a flat surface with a pair of singularities of order $k_1$ and $k_2$ and a vertical saddle connection of length $\varepsilon$. We will denote be $\Psi(S,\gamma,\varepsilon)$ this surface.
The construction is continuous with respect to the variations of $\varepsilon$.

\subsection{Dependence on small variations of the path}
The previous construction might depend on the choice of the broken line.  We show the following proposition:

\begin{prop}  \label{indep_petite_variation}
 Let $\gamma$ and $\gamma^{\prime}$ be two broken lines that both start from $P$, sector $I$ and end to $P$, sector $II$. Let $\varepsilon$ be a positive real number. We assume that there exists an open subset $U$ of $S$, such that:
 \begin{itemize}
 \item $U$ contains  $\gamma\backslash \{P\}$ and $\gamma^{\prime}\backslash \{P\}$. 
\item $U$ is homeomorphic to a disc and have no conical singularities.
\item The surgery described in section \ref{break_non_loc}, with parameters $(\gamma,\varepsilon)$ or $(\gamma^{\prime},\varepsilon)$ does not affect $\partial U \backslash P$.
\end{itemize}
Then $\Psi(S,\gamma,\varepsilon)$ and $\Psi(S,\gamma^{\prime},\varepsilon)$ are isometric.
 \end{prop}

 \begin{figure}[htbp]
 \begin{center}
 \input{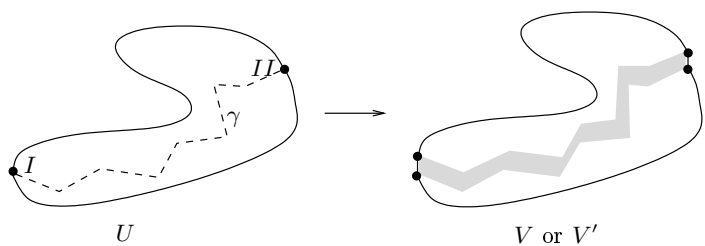}
 \caption{The boundary of $U$ and  $V$ (or  $V^{\prime}$).}
\end{center}
\end{figure}

\begin{proof} 
We denote by $\partial U$ the boundary of the natural compactification of $U$ (that differ from the closure of $U$ in $S$, see section \ref{homolo}). We denote by $\tilde {P}$ and $\tilde{P}^{\prime}$ the ends of $\gamma$ in $\partial U$  (that are also the ends of $\gamma^{\prime}$ by assumption). We denote by $V$ (\emph {resp.}  $V^{\prime}$) the flat discs obtained from $U$ after the hole surgery along $\gamma$ (\emph{resp.} $\gamma^{\prime}$). Our goal is to prove that $V$ and $V^{\prime}$ are isometric.

The hole surgery along $\gamma$ (\emph{resp.} $\gamma^{\prime}$) does not change the metric in a neighborhood of $\partial U \backslash \{\tilde{P},\tilde{P}^{\prime}\}$. Furthermore, the fact that both $\gamma$ and $\gamma^{\prime}$ starts and ends at sectors $I$ and $II$ correspondingly implies that $V$ and $V^{\prime}$ are isometric in a neighborhood of their boundary. We denote by $f$ this isometry. Surprisingly, we can find two flat discs that are isometric in a neighborhood of their boundary but not globally isometric  (see Figure \ref{disques_plats}).

\begin{figure}[htpb]
\begin{center}
\includegraphics[scale=0.5]{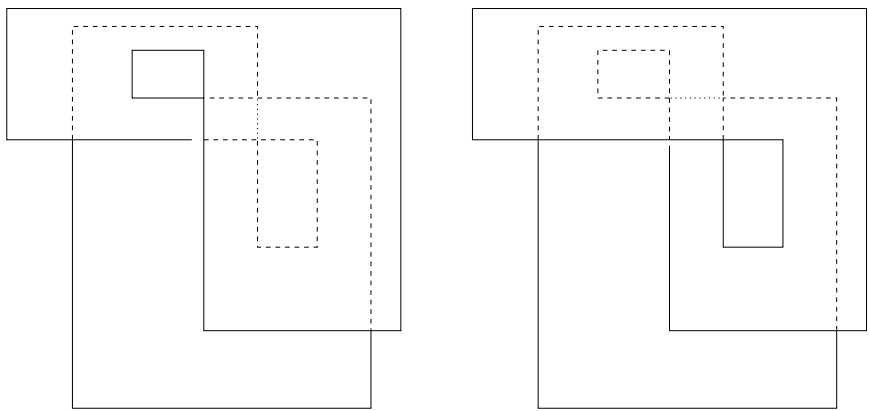}
\caption{Immersion in $\mathbb{R}^2$ of two non isometric flat discs with isometric boundaries.}
\label{disques_plats}
\end{center}
\end{figure}

In our case, we have an additional piece of information that will make the proof possible: hole transport does not change the vertical foliation  (recall that the hole is always assumed to be vertical). Therefore, for each vertical geodesics in $V$ with end points $\{x, y\}\subset \partial V$, then $\{f(x),f(y)\}$ are the end points of a vertical geodesic of $V^{\prime}$.

 \begin{figure}[htbp]
 \begin{center}
 \input{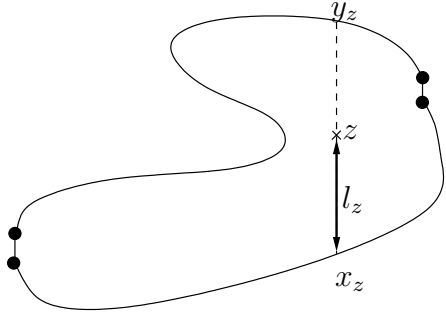}
 \caption{Parameters on a flat disc.}
 \label{parametre:disques}
\end{center}
\end{figure}

For each $z\in V$ we define $x_z\in \partial V$ (\emph{resp.} $y_z$)  the intersection of the vertical geodesic starting from $z$ in the negative direction (\emph{resp.} positive direction) and the boundary of $V$ (see Figure \ref{parametre:disques}). We also call $l_z$ the length of this geodesic. We can assume that $\partial V$ is piecewise smooth. So we can restrict ourself to the open dense subset $V_1 \subset V$  of $z$ such that $x_z$ and $y_z$ are regular and nonvertical points. 

Then we define $\Phi: V_1 \rightarrow V^{\prime}$ that send $z$ to $\phi_{l_z}(f(x_z))$, where, $\phi$ is the vertical geodesic flow. Because $V$and $V^{\prime}$ are translation structures, the length of the vertical segment $[x_z,y_z]$ is obtained by integrating the corresponding $1-$form along any path between $x_z$ and $y_z$. Such a path can be chosen in a neighborhood of the boundary of $V$. Then, the isometry $f$ implies that this length is the same as the length of the vertical segment $[f(x_z),f(y_z)]$. Therefore $\Phi$ is well defined and coincide to $f$ in a neighborhood of the boundary of $V$. This map is also smooth because $z\mapsto (x_z, l_z)$ are smooth on $V_1$. It's easy to check that $D\Phi(z)\equiv Id$ and that $\Phi$ continuously extends to an isometry from $V$ to $V^{\prime}$.
\end{proof}


\begin{cor}\label{petites:variations}
Let $\gamma^{\prime}$ be close enough to $\gamma$ and such that $\gamma$ and $\gamma^{\prime}$ intersect the same sectors of a neighborhood of $P$. Then $\Psi(S,\gamma,\varepsilon)$ and $\Psi(S,\gamma^{\prime},\varepsilon)$ are isomorphic for $\varepsilon$ small enough. 
\end{cor}

\begin{proof}
If $\gamma^{\prime}$ is close enough to $\gamma$ (and intersect the same sectors in a neighborhood of $P$), then there exists a open flat disk that contains $\gamma$ and $\gamma^{\prime}$.
\end{proof}

\begin{rem}
Using proposition \ref{indep_petite_variation}, one can also extend hole transport along a differentiable curve.
\end{rem}


\section{Simple configuration domains} \label{simple_cusps}

We recall the following notation: if $\mathcal{Q}(k_1,k_2,\dots,k_r)$  is a stratum of meromorphic quadratic differentials with at most simple poles, then $\mathcal{Q}_1(k_1,k_2,\dots,k_r)$ is the subset of area 1 flat surfaces in $\mathcal{Q}(k_1,k_2,\dots,k_r)$, and $\mathcal{Q}_{1,\delta}(k_1,k_2,\dots,k_r)$  is the subset of flat surfaces in $\mathcal{Q}_1(k_1,k_2,\dots,k_r)$ that have at least  a saddle connection of length less than $\delta$.

The goal of this section is to prove the following theorem: 
\begin{ths}
\label{cusp_elementaire} 
Let $\mathcal{Q}(k_1,k_2,\dots,k_r)$ be a stratum of quadratic differentials with $(k_1,k_2)\neq(-1,-1)$ and such that the stratum $\mathcal{Q} (k_1+k_2,k_3, \dots,k_r)$ is connected.  Let $\mathcal{C}$ be the subset of flat surfaces $S$ in $\mathcal{Q}^{N} (k_1,\dots,k_r)$ such that the shortest saddle connection of $S$ is simple and joins a singularity of order $k_1$ to a singularity of order $k_2$. For any pair $N\geq 1$ and $\delta>0$, the sets $\mathcal{C}$, $\mathcal{C} \cap \mathcal{Q}_1(k_1,k_2,\dots,k_r)$ and $\mathcal{C} \cap \mathcal{Q}_{1,\delta}(k_1,k_2,\dots,k_r)$ are non empty and connected. 
 
\end{ths}

In this part we denote by $P_1$ and $P_2$ the two zeros of order $k_1$ and $k_2$ respectively and by $\gamma$ the simple saddle connection between them. There are two different cases.
\begin{itemize} 
\item When $k_1$ and $k_2$ are not both odd, then
 there exists a canonical way of shrinking the saddle connection $\gamma$ if it is small enough. Furthermore,
 this surgery doesn't change the metric outside a neighborhood of $\gamma$. This
 is the local case.  
\item When $k_1$ and $k_2$ are both odd, then we still can shrink
 $\gamma$, to get a surface in the stratum  $\mathcal{Q}(k_1+k_2,k_3,\ldots,k_r)$, but this changes the metric  outside a neighborhood of $\gamma$ and this is not canonical. This is  done by reversing the procedure of section \ref{break_non_loc}. 
\end{itemize}

\subsection{Local case} 
\subsubsection{Breaking up a singularity} \label{break_loc}
Here we follow \cite{EMZ,MZ}. Consider a singularity $P$ of order $k\geq 0$, and a partition $k=k_1+k_2$ with $k_1,k_2\geq -1$. We assume that $k_1$ and $k_2$ are not both odd. If $\rho$ is small enough, then the set $\{x\in S, d(x,P)<\rho \} $ is a metric disc embedded in $S$.  It is obtained by gluing  $k+2$ standards half-disks of radius $\rho$.

\begin{figure}[htbp]
%
\begin{center}
\input{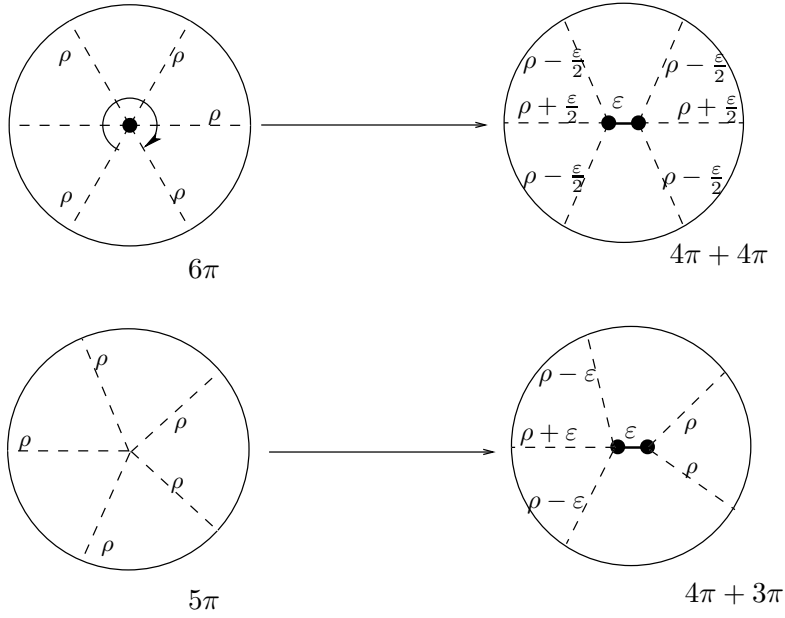}
\caption{Breaking up a zero into two zeroes (after~\cite{EMZ,MZ}).}
\label{break_sing_loc}
\end{center}
\end{figure}

There is a well known local construction  that breaks the singularity $P$ into two singularities of order $k_1$ and $k_2$, and which is obtained by changing continuously the way of gluing the half-discs together (see figure \ref{break_sing_loc}, or \cite{EMZ,MZ}). This construction is area preserving.

\subsubsection{Structure of the neighborhood of the principal boundary} 
When $\gamma$ is small enough, (for example $|\gamma|\leq |\gamma^ {\prime}|/10$, for any other saddle connection $\gamma^{\prime}$), then we can perform the reverse construction because a neighborhood of $\gamma$ is precisely obtained from a collection of half-discs glued as before. This defines a canonical map $\Phi:V\rightarrow \mathcal{Q}(k_1+k_2,k_3, \ldots,k_r)$, where $V$ is a subset of $\mathcal{Q}(k_1,k_2,k_3, \ldots,k_r)$. We can choose $U^N\subset V$ such that $\Phi^{-1}(\{\widetilde{S}\})\cap U^N$ is the set of surfaces such that the shrinking process leads to $\widetilde{S}$, and whose smallest saddle connection is of length smaller than $\min(\frac{\tilde{|\gamma|}}{100}  ,  \frac{\tilde{|\gamma|}}{2N}))$ with $\tilde{\gamma}$ the smallest saddle connection of $\widetilde{S}$.  From the proof of Lemma $8.1$ of \cite{EMZ}, this map gives to $U^N$ a structure of a topological orbifold bundle over $\mathcal{Q}(k_1+k_2,k_3,\ldots,k_r)$, with the punctured disc as a fiber.
By assumption, $\mathcal{Q}(k_1+k_2,k_3,\ldots,k_r)$ is connected, and therefore $U^N$ is connected, so the proof will be completed after the following three steps:
\begin{itemize}
\item $U^N\subset \mathcal{C}$. 
\item There exists $L>0$ such that $Q^L(k_1,\dots,k_r)\cap \mathcal{C} \subset U^N$.
\item For any $S\in \mathcal{C}$, there exists a continous path $(S_t)_t$ in $\mathcal{C}$ that joins $S$ to $Q^L(k_1,\dots,k_r)$.
\end{itemize}

\subsubsection{Proof of Theorem \ref{cusp_elementaire}: local case}
To prove the first step, it is enough to show that $U^N$ is a subset of $\mathcal{Q}^N(k_1,k_2,\dots,k_r)$:  let $S$ be a flat surface in $U^N$ and let $\widetilde{S}=\Phi(S)$. We denote by $\gamma$ the smallest saddle connection of $S$. The surgery doesn't change the surface outside a small neighborhood of the corresponding singularity of  $\widetilde{S}$. If $|\tilde{\gamma}|$ is the length of the smallest saddle connection of $\widetilde{S}$, then $S$ has no saddle connection of length smaller than $\tilde{|\gamma|}-|\gamma|$ except $\gamma$, which has length smaller than $\frac{\tilde{|\gamma|}}{2N}$ by construction. 
We have 
$\frac{|\tilde{\gamma}|-|\gamma|}{|\gamma |}=\frac{|\tilde{\gamma}|}{|\gamma|}-1>2N-1\geq N$, so $S$ belongs to $\mathcal{Q}^N(k_1,k_2,\dots,k_r)$.
Hence we have proved that $U^N\subset \mathcal{C}$.

To prove the second step, we remark that if $S\in Q^L(k_1,\dots,k_r)\cap \mathcal{C}$, for $L\geq 10$, then the smallest saddle connection of $\Phi(S)$ is of length at least $L|\gamma|-|\gamma|$, where $\gamma$ is the smallest saddle connection of $S$. Hence if $|\gamma|\leq \min(\frac{(L-1)|\gamma|}{100}, \frac{(L-1)|\gamma|}{2N})$ then $S\in U^N$. So we have proved that $Q^L(k_1,\dots,k_r)\cap \mathcal{C}\subset U^N$ for $L\geq \max(101,2N+1)$.

The last step is given by the following lemma:
\begin{lem} 
Let $S$ be a surface in $\mathcal{Q}^N  (k_1,\dots,k_r)$ whose
smallest saddle connection $S$ is simple and joins a singularity of order $k_1$ to a singularity of order $k_2$,  and let $L$ be a positive number. Then we can find a continuous path in $\mathcal{Q}^N  (k_1,\dots,k_r)$, that joins $S$ to a surface whose second smallest saddle connection is at least $L$ times greater than the smallest one.
\end{lem}

\begin{proof} 
 The set $\mathcal{Q}^N(k_1,\dots,k_r)$ is open, so up to a small continuous perturbation of $S$, we can  assume that $S$ has no vertical saddle connection except the smallest one.

Now we apply the geodesic flow $g_t$ to $S$. There is a natural bijection from the saddle connections of $S$ to the saddle connections of $g_t. S$. The holonomy vector $v=(v_1,v_2)$ of a saddle connection  becomes $v_t=(e^{-t} v_1,e^t v_2)$. This imply that the quotient of the length of a given saddle connection to the length of the smallest one increases and goes to infinity.

The set of holonomy vectors of saddle connections is discrete, and therefore, any other saddle connection of $g_t.S$ has length greater than $L$ times the length of the smallest one, as soon as $t$ is large enough.
\end{proof}

Note that the previous proof is the same if we restrict ourselves to area 1 surfaces. The case for when restricted to the $\delta$-neighborhood of the boundary is also analogous, since $U^N \cap \mathcal{Q}_{1,\delta}(k_1,\dots,k_r)$ is still a bundle over $\mathcal{Q}_{1}(k_1,\dots,k_r)$ with the punctured disc as a fiber.

Hence the theorem is proven when  $k_1$ and $k_2$ are non both odd.

\subsection{Proof of theorem \ref{cusp_elementaire}: non-local case}\label{non_loc}
We first show that two surfaces that are close enough to the stratum
$\mathcal{Q}(k_1+k_2,k_3,\ldots,k_r)$ (in a certain sense that will be specified below) belong to the same configuration domain. Then we show that we can always continuously reach that neighborhood.

\subsubsection{Neighborhood of the principal boundary} 

Contrary to the local case, we do not have a canonical map from a subset of $\mathcal{Q}(k_1,k_2,\ldots,k_r)$  to $\mathcal{Q}(k_1+k_2,\ldots,k_r)$ that gives to this subset a structure of a bundle.

Let $S \in \mathcal{Q}(k_1+k_2,\ldots,k_r)$, and let $\nu$ be a path in $S$, we will say that $\nu$ is \emph{admissible} if it satisfies the hypothesis of the singularity breaking procedure of section \ref{break_non_loc}. Let $\nu$ be an admissible closed path whose end point is a singularity $P$ of degree $k_1+k_2$  and let $\varepsilon>0$ be small enough for the breaking  procedure. Recall that  $\Psi(S,\nu,\varepsilon)$ denotes the surface in $\mathcal{Q}(k_1,k_2,\ldots,k_r)$ obtained  after breaking the singularity $P$, using the procedure of section \ref{break_non_loc} along the path $\nu$, with a vertical hole of length $\varepsilon$. 

\begin{prop} \label{relev_chemin}
 Let $(S,S^{\prime})$ be a pair of surfaces in $\mathcal{Q}(k_1+k_2,\ldots,k_r)$ and $\nu$ (\emph{resp.} $\nu^{\prime}$) be an admissible  broken line in $S$ (\emph{resp.} $S^{\prime}$). Then
$\Psi(S,\gamma,\varepsilon)$ and $\Psi(S^{\prime},\gamma^{\prime},\varepsilon)$ belong to the same configuration domain for any sufficiently small  $\varepsilon$.  
\end{prop}

\begin{proof}
By assumption, $\mathcal{Q}(k_1+k_2,\ldots,k_r)$ is connected, so
there exists a path $(S_t)_{t\in[0,1]}$, that joins $S$ and $S^{\prime}$. We can find a family of broken lines  $\gamma_t$ of $S_t$ such that, for $\varepsilon$ small enough, the map $t\mapsto \Psi(S_t,\gamma_t,\varepsilon)$ is well defined and continuous for $t\in[0,1]$. The surface $\Psi(S^{\prime},\gamma_1,\varepsilon)$ might differ from
$\Psi(S^{\prime},\gamma^{\prime},\varepsilon)$ for two reasons: 
\begin{itemize}
 \item The paths $\gamma_1$ and $\gamma^{\prime}$, that both start from the same singularity $P$ , might not start and end at the same sectors. In that case, we consider the path $r_{\theta} S^{\prime}$ obtained by rotating the surface $S^{\prime}$ by an angle of $\theta$. We find as before a family of broken lines $\gamma_{1,\theta} \in r_{\theta} S^{\prime}$. Then, for some  $\theta_k$ an integer multiple of $\pi$,  we will have $r_{\theta_k} S^{\prime}=S^{\prime}$ and $\gamma_{1,\theta_k}$ that starts and ends on the same sectors than $\gamma^{\prime}$.  
\item Even if the paths $\gamma_1$ and $\gamma^{\prime}$ start and end in the same sectors of the singularity $P$, they might be very different (for example in a different homotopy class of $S^{\prime}\backslash \{sing\}$), so  Proposition \ref{indep_petite_variation} does not apply. This case is solved by the following lemma, which says that the resulting surfaces are in the same configuration domain.
\end{itemize} 
\end{proof}

\begin{lem} \label{independant} 
For any surface $S\in \mathcal{Q}(k_1+k_2,k_3,\dots,k_r)$, the configuration domain that contains a surface obtained by the non-local singularity breaking construction does not depend on the choice of the admissible path, once sector $I$ is chosen, and the hole is small enough.  
\end{lem}

\begin{proof} 
We consider a surface $S$ in $\mathcal{Q}(k_1+k_2,\dots,k_r)$ and perform the breaking procedure. We do not change the resulting configuration domain if we perform some small perturbation of $S$. Therefore, we can assume that  the vertical flow of $S$ is minimal (this is the case for almost all surface).  Now we consider an admissible path and perform the corresponding singularity breaking procedure and get a surface $S_1$.  Then we consider a horizontal segment in sector $I$ adjacent to the singularity $k_1$.  Then we perform the same construction for another admissible path (and get a surface $S_2$) and consider a horizontal segment of the same length as before.

Because the hole transport preserves the vertical foliation,  the first return maps of the vertical flow in the two surfaces are the same as soon as the hole is small enough.

Now from Lemma \ref{key_fam}, there exists $(\sigma,l)$ such that $S_1$ and $S_2$ belong to $ZR^{\prime}(Q_{\sigma,l}^{\prime})$, with parameters $\zeta^1_1,\ldots,\zeta^1_{l+m}$ and $\zeta^2_1,\ldots,\zeta^2_{l+m}$. Note that $Re(\zeta^1_i)=Re(\zeta^2_i)$, because these depends only on the first returns maps of the vertical flow (and they coincide). The family of polygons with parameters $t \zeta^1_i+(1-t)\zeta^2_i$ gives a path in $MZ^{\prime}(Q_{\sigma,l}^{\prime})$ that joins $S_1$ and $S_2$. Furthermore, the singularity breaking procedure is continuous with respect to $\varepsilon$. Hence,  for all $i$, $\zeta^1_i$ and $\zeta^2_i$ are arbitrary close as soon as $\varepsilon$ is small enough. Consequently, the constructed  path in $MZ^{\prime}(Q_{\sigma,l}^{\prime})$  keeps being in a configuration domain.

\begin{figure}
\begin{center}
\input{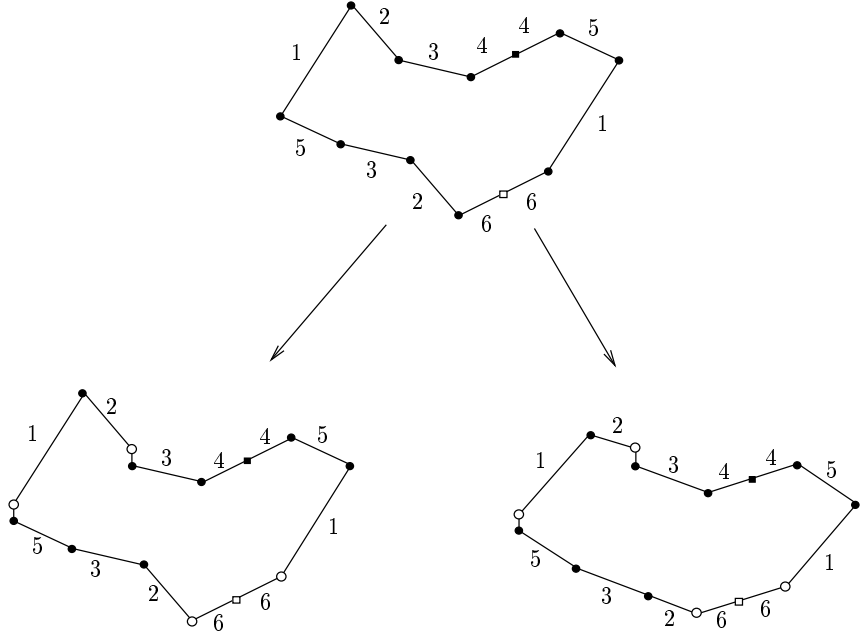}
\caption{Breaking a singularity with two different paths.}
\end{center}
\end{figure}

\end{proof}

Now for each $S\in Q(k_1+k_2,\dots,k_r)$ and each admissible path
$\gamma$, we can find $\varepsilon_{S,\gamma}$ maximal  such that
$\Psi(S,\gamma,\varepsilon)\in Q^N(k_1,\dots,k_r)$ for all
$\varepsilon < \varepsilon_{S,\gamma}$. Now we consider the set 
\[U^N=\bigcup_{\theta \in [0,2\pi]} \bigcup_{S,\gamma}\ 
\bigcup_{0<\varepsilon< \varepsilon_{S,\gamma}}
r_{\theta}(\Psi(S,\gamma,\varepsilon))\] 
This subset of $Q^N(k_1,\dots,k_r)$ is a subset of a configuration domain from Proposition \ref{relev_chemin}.

\subsubsection{Reaching a neighborhood of the principal boundary}

Now we consider a surface in $\mathcal{Q}^N(k_1,\ldots,k_r)$ whose unique smallest saddle connection joins a singularity of order $k_1$ to a singularity of order $k_2$. As in the local case, we can assume that its smallest saddle connection is vertical and that there are no other vertical saddle connections. Then we apply the Teichmüller geodesic flow. This allows us to assume that the smallest saddle connection is arbitrary small compared to any other saddle connection. 

We then want to contract the saddle connection using the reverse
procedure of section \ref{break_non_loc}. 

\begin{prop}\label{contractN} 
Let $N$ be greater than or equal to $1$. There exists $L>N$ such that 
$Q^L(k_1,\dots,k_r)\cap \mathcal{C}\subset U^N$.
\end{prop}

\begin{proof} 


We choose $L$ large enough such that we can find $L^{\prime}$ satisfying  $2N<L^{\prime}$, and $1\ll L^{\prime}\ll L$.  Denote by $\gamma$ the smallest saddle connection and by  $\varepsilon$ its length. We want to find a path suitable for reversing the construction of section \ref{break_non_loc}. When contracting $\gamma$ in such way, we must insure that the surface stay in $Q^N(k_1,\dots,k_r)$, by keeping a lower bound of the length of the saddle connections different from the shortest one.

Let $B$ be the open $L^{\prime}\varepsilon$-neighborhood of $\gamma$, and $\{B_i\}_{i\in \{3,\dots,r\}}$ the open $L^{\prime}\varepsilon$-neighborhoods of the singularities that are not end points of $\gamma$. Note that each of these neighborhoods is naturally isometric to a collection of half-disk glued along their boundary. We denote by  $S^{\prime}$ the closed subset of $S$ obtained by removing to $S$ the set $\cup_i B_i \cup B$. 

\begin{figure}[htb]
\begin{center}
\input{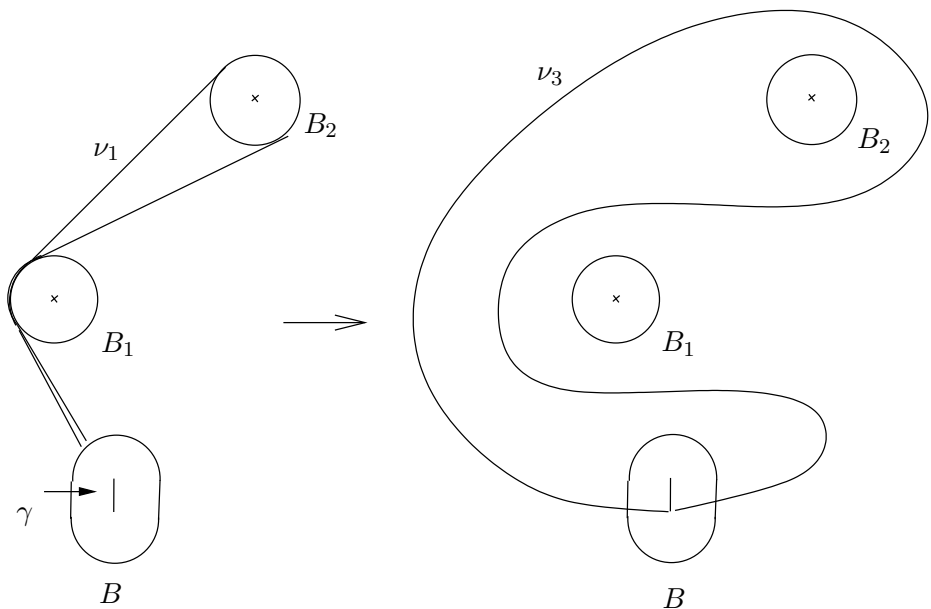}
\caption{Constructing a suitable path.}
\end{center}
\end{figure}

Now we consider the set of paths of $S^{\prime}$ whose end points are on  $\partial B$ and with nontrivial holonomy (which makes senses in a  neighborhood of $\partial B$), and we choose a path $\nu_1$ of minimal length with this property. Note that, we do not change the holonomy of a path by ``uncrossing'' generic self intersections (see figure \ref{uncross}). Therefore, we can choose our path such that, after a small perturbation, it has no self intersections.

 \begin{figure}[htb]
\begin{center}
\input{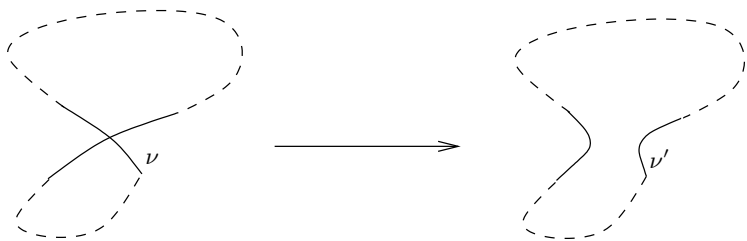}
\caption{Uncrossing an intersection does not change the holonomy.}
\label{uncross}
\end{center}
\end{figure}

 Now the condition $L^{\prime}\ll L$ implies that we can find find a path $\nu_2$ in the same homotopy class, such that the $\varepsilon$-neighborhood of $\nu_2$ is homeomorphic to a disk. Now joining carefully the end points of $\nu_2$ to each sides of $\gamma$, we get a path $\nu_3$. By construction, we can use this path to contract the saddle connection $\gamma$. The surgery doesn't touch the $\varepsilon N$-neigborhoods of the singularities, except
 for the end points of $\gamma$, hence any saddle connection that starts from such singularity will have a length greater than $N\varepsilon$ during the shrinking process. 
 A saddle connection starting from an end point a $\gamma$, and different from $\gamma$ will leave~$B$. Choosing properly $\nu_3$, then the length of such  saddle connection will have a length greater than $(L^{\prime}-1)\varepsilon$ 
 during the shrinking process, and $L^{\prime}-1\geq N+(N-1)\geq N$.

Therefore, when contracting $\gamma$, there is no saddle connection except $\gamma$ that is of length smaller than $N|\gamma| \leq N\varepsilon$, were $\varepsilon$ is the initial length of the  saddle connection $\gamma$.  Up to rescaling the surface, we can assume that the area of the surface is constant under the deformation process.
\end{proof}

Now let $\mathcal{C}$ be the open subset of surfaces in  $\mathcal{Q}^N(k_1,\ldots,k_r)$ whose unique smallest saddle connection joins a singularity of order $k_1$ to a singularity of order $k_2$. The previous proposition shows that there exists a path from any $S\in \mathcal{C}$ to $U^N$, which is pathwise connected. Therefore $\mathcal{C}$ is pathwise connected and hence, connected. Then we have proven the theorem for the case when $k_1$ and $k_2$ are odd.


\section{Configuration domains in strata of quadratics differentials on the Riemann sphere} \label{CP1}

\begin{figure}
\begin{center}
\input{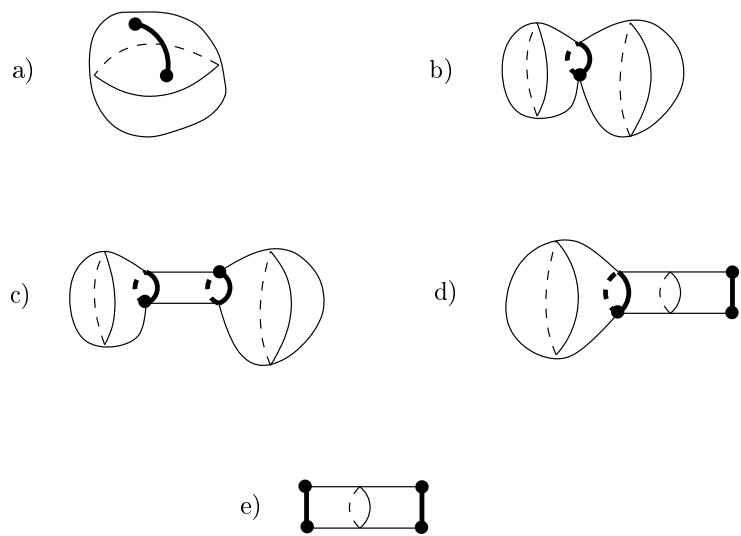}
\caption{``Topological picture'' of  configurations for $\mathbb{CP}^1$.}
\label{picture:cp1}
\end{center}
\end{figure}

In \cite{B} we proved Theorem \ref{th:config} describing all the  configurations of \^homologous saddle connections that exist on a given stratum of quadratic differential on $\mathbb{CP}^1$. We now show that they are in bijections with the configuration domains. In this section, we denote by $\gamma$ a collection $\{\gamma_i\}$ of saddle connections.
\begin{ths}\label{th:config}
Let $\mathcal{Q}(k_1^{\alpha_1},\ldots,k_r^{\alpha_r}, -1^s)$ be a stratum of quadratic differentials on $\mathbb{CP}^1$ different from $\mathcal{Q}(-1^4)$, and let $\gamma$ be a maximal collection of \^homologous saddle connections on a generic surface in that stratum. Then the possible configurations for $\gamma$ are given in the list below (see Figure \ref{picture:cp1}).
\begin{itemize}
\item[a)]
Let $\{k,k^{\prime}\}\subset \{k_1^{\alpha_1},\ldots,k_r^{\alpha_r}, -1^s\}$ be  an unordered pair of integers with $(k,k^{\prime})\neq (-1,-1)$. The set $\gamma$ consists of a single saddle connection joining a singularity of order
$k$ to a distinct singularity of order $k^{\prime}$.
\item[b)]
Consider $(a_1,a_2)$ a pair of strictly positive integers such that $a_1+a_2=k\in\{k_1, \ldots,k_r\}$ (with $k\neq 1$), and a partition $A_1 \sqcup A_2$ of $\{k_1^
{\alpha_1},\ldots,k_r^{\alpha_r}\} \backslash \{k\}$. The set $\gamma$ consists of a simple saddle connection that decomposes the sphere into two 1-holed spheres $S_1$ and $S_2$, such that each $S_i$ has interior singularities of positive order given by $A_i$ and $s_i=(\sum_{a\in A_i}a)+a_i+2$ poles, and has a~single boundary singularity of order $a_i$.
\item[c)]
Consider $\{a_1,a_2\}\subset \{k_1^{\alpha_1},\ldots,k_r^{\alpha_r}\}$  a pair of integers. Let $A_1 \sqcup A_2$ be a partition of $\{k_1^{\alpha_1},\ldots,k_r^{\alpha_r}\}\backslash \{a_1,a_2\}$. The set $\gamma$ consists of two closed saddle connections that decompose the sphere into two 1-holed spheres $S_1$ and $S_2$ and a  cylinder,  and such that each $S_i$ has interior singularities of positive orders given by $A_i$ and  $s_i=(\sum_{a\in A_i}a)+a_i+2$ poles and has a boundary singularity of order $a_i$.
\item[d)]
Let $k\in \{k_1,\ldots,k_r\}$.  The set $\gamma$ is a pair of saddle connections of different lengths, and such that the largest one  starts and ends from a singularity of order $k$ and  decompose the surface into a 1-
holed sphere and a half-pillowcase, while the shortest one joins a pair of poles and lies on the other end of the half pillowcase.
\end{itemize}
When the stratum is $\mathcal{Q}(-1^4)$, there is only one configuration, which correspond to two saddle connections are the two boundary components of a cylinder (the surface is a ``pillowcase'', see Figure~\ref{exemple_trivial}).
\end{ths}

The following lemma (due to Kontsevich) implies that Theorem \ref{cusp_elementaire} will apply for any stratum of quadratic differentials on $\mathbb{CP}^1$ (see also~\cite{KoZo}).

\begin{NNlem}[Kontsevich]
 Any stratum of quadratic differentials on $\mathbb{CP}^1$ is non
 empty and connected.  
\end{NNlem}
\begin{proof} 
 There is only one complex structure on $\mathbb{CP}^1$.  Therefore,
 we can work on the standard atlas  $\mathbb{C} \cup \left (\mathbb{C}^*\cup \infty\right) $ of the Riemann  sphere.

Now we remark that if we fix $(z_1,\dots,z_r)\in\mathbb{C}^r$ that are pairwise distinct, and $k_1,\dots,k_r$ some integers greater than or equal to $-1$, then the quadratic differential on $\mathbb{C}$,  $q(z)=\prod(z-z_i)^{k_i} dz^2$, extends to a quadratic differential on $\mathbb{CP}^1$ with possibly a singularity of order $-4-\sum_i k_i$ over the point $\infty$. Now two quadratic differentials on a compact Riemann surface with the same singularity points are equal up to a multiplicative constant (because they differ by a holomorphic function). 

 Therefore, any strata of quadratic differentials on $\mathbb{CP}^1$ is a quotient of  $\mathbb{C}$ times a space of configurations of points on a sphere,  which is connected.
\end{proof}

Now let $S \in \mathcal{Q}^N(k_1^{\alpha_1},\ldots,k_r^{\alpha_r}, -1^s)$ . We can define  $\mathcal{F}_S$ to be the maximal collection of \^homologous saddle connections that contains the smallest one. We have the following lemma.

\begin{lem}\label{lem_config}
The configuration associated to $\mathcal{F}_S$ is locally constant with respect to $S$. 
\end{lem}

\begin{figure}[ht]
\begin{center}
\input{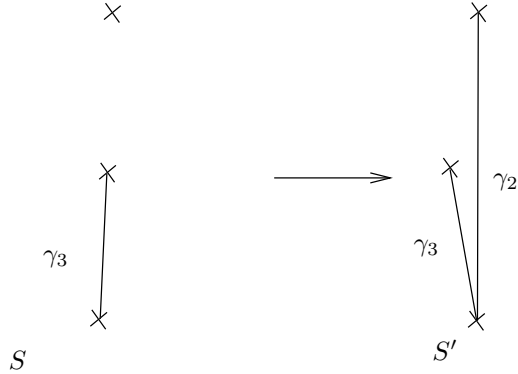}
\caption{The configuration associated to $\mathcal{F}_S$ is locally constant.}
\label{loc_cst}
\end{center}
\end{figure}

\begin{proof}
Any saddle connection in $\mathcal{F}_S$ persists under small deformation. This lemma is obvious as soon the number of elements of  $\mathcal{F}_S$ is locally constant.

Let $\gamma_1$ be a saddle connection of minimal length. We assume that after a small perturbation $S^{\prime}$ of $S$, we get a bigger  collection of saddle connections. That means that a new saddle connection $\gamma_2$ appears. Therefore there was another saddle connection $\gamma_3$ non\^homologous to $\gamma_1$, of length less than or equal to $|\gamma_2/2|$ (see figure \ref{loc_cst}). But this is impossible since it would therefore be of length less than or equal to the length of $\gamma_1$, contradicting the hypothesis.
\end{proof}

\begin{MT} 
Let $\mathcal{Q}(k_1^{\alpha_1},\ldots,k_r^{\alpha_r}, -1^s)$ be a stratum of quadratic differentials with at most simple poles.  Let $N$ be greater than or equal to~$1$. There is a natural bijection between the configurations of \^homologous saddle connections on $\mathcal{Q}(k_1^{\alpha_1},\ldots,k_r^{\alpha_r}, -1^s)$ described in Theorem \ref{th:config} and the connected components of $\mathcal{Q}^N(k_1^{\alpha_1},\ldots,k_r^{\alpha_r}, -1^s)$ 
\end{MT}

\begin{proof} 
Lemma \ref{lem_config} implies that there is a well defined map $\Psi$ from the set of connected components of $\mathcal{Q}^N(k_1^{\alpha_1},\ldots,k_r^{\alpha_r}, -1^s)$ to the set of existing configurations for the stratum. This map is surjective because if we choose a generic surface $S$ with a maximal collection of \^homologous saddle connections $\gamma$ that realizes the given configuration $\mathcal{C}$, then after a small continuous perturbation of the surface, we can assume that there are no other saddle connections on $S$ parallel to an element of $\gamma$. Then we apply the Teichmüller geodesic flow to contract the elements of  $\gamma$, until $\gamma$ contains the smallest saddle connection of the surface. Then by construction, this surface belongs to $\Psi^{-1}(\mathcal{C})$.

Now we prove that $\Psi$ is injective.  
We keep the notations of the previous theorem, and consider $U=\Psi^{-1}(\{\mathcal{C}\})$, for  $\mathcal{C}$ any existing configuration:
 
-If $\mathcal{C}$ belongs to the a) case, then $U$ is connected from  Theorem \ref{cusp_elementaire} and the lemma of Kontsevich.

-If $\mathcal{C}$ belongs to the b) case, then  we consider a surface $S$ in $U$. Its smallest saddle connection  $\gamma_0$ is closed and separates the surface in a pair $(S_1,S_2)$ of 1-holed spheres with  boundary  singularities of orders $a_1$ and $a_2$ correspondingly. Now for each $S_i$ we decompose the boundary saddle connection of $S_i$ in two segments starting from the boundary singularity, and glue together these two segments, then 
we get a pair of closed flat spheres  $\widehat{S}_i\in \mathcal{Q}(A_i,a_i-1,-1^{s_i+1})$, $i=1,2$. For each of the sphere, the smallest saddle connection $\gamma_i^{\prime}$ is simple and  joins a singularity $Q_i$  of order $(a_i-1)$ to a newborn pole $P_i$, and is of length $|\gamma_0|/2$, where $|\gamma_0|$ is the length of $\gamma_0$. Let $\eta_i$ be the smallest saddle connection of $\widehat{S}_i$ except $\gamma_i^{\prime}$.
\begin{itemize}
\item If $\eta_i$ intersects the interior of $\gamma_i^{\prime}$, then it is easy to find another saddle connection on $\widehat{S}_i$, smaller than $\eta_i$ and different from $\gamma_i^{\prime}$. 
\item If $\eta_i$ does not intersect $\gamma_i^{\prime}$, or intersect it in $Q_i$, then $\eta_i$ was a saddle connection on $S$, hence $|\eta_i|> 2N |\gamma_i^{\prime}|$.
\item If $\eta_i$ intersects $P_i$, then we can find a saddle connection in $S$ of length smaller than $|\eta_i|+|\gamma_0|/2$.
\end{itemize}

These remarks imply hat $\widehat{S}_i$ is in  $\mathcal{Q}^{2N-1}(A_i,a_i-1,-1^{s_i+1})$ which is a subset of $\mathcal{Q}^{N}(A_i,a_i-1,-1^{s_i+1})$. Hence we have defined a map $f$ from $U$ to $U_1\times U_2$, with $U_i$ a simple configuration domain of $Q^N(A_i,a_i-1,-1^{s_i+1})$. 

Conversely,  let $\{\widehat{S}_i\}_{i\in \{1,2\}} $ be two surfaces in $\mathcal{Q}^{2N}(A_i,a_i-1,-1^{s_i+1})$, such that for each $\widehat{S}_i$, the smallest saddle connection $\gamma_i$ is simple and joins a pole to a singularity of order $a_i-1$. If $\gamma_1$ and $\gamma_2$ are in the same direction and have  the same length, then we can reconstruct a surface $S=f^{-1}(\widehat{S}_1,\widehat{S}_2)$ in $\mathcal{Q}(k_1^{\alpha_1},\ldots,k_r^{\alpha_r}, -1^s)$  by cutting $\widehat{S}_i$ along $\gamma_i$, and gluing together the two resulting surfaces. Here, $S$  belongs to $\mathcal{Q}^N(k_1^{\alpha_1},\ldots,k_r^{\alpha_r}, -1^s)$. 
Note that in the reconstruction of the surface, the length of smallest saddle connection is doubled, hence we must start from $\mathcal{Q}^{2N}(A_i,a_i-1,-1^{s_i+1})$, and not $\mathcal{Q}^{N}(A_i,a_i-1,-1^{s_i+1})$.

Now we proove the connectedness of $U$: let $X^1, X^2$ be two flat surfaces in $U$. After a small perturbation and after applying the geodesic flow, we get a surface $S^1$ (\emph{resp.} $S^{2}$) in the same connected component of $U$ as $X^1$ (\emph{resp} $X^{2}$), with $S^1$ and $S^{2}$ in $\mathcal{Q}^{2N}(k_1^{\alpha_1},\ldots,k_r^{\alpha_r}, -1^s)$. 

There exists  continuous paths $(S_{i,t})_{t\in[1,2]}\in \mathcal{Q}^{2N}(A_i,a_i-1,-1^{s_i+1})$ such that $(S_{1,j},S_{2,j})=f(S^j)$ for $j=1,2$. The pair $(S_{1,t},S_{2,t})$ belongs to $f(U)$ if and only if their smallest saddle connections are parallel and have the same length. 
This condition is not necessary satisfied, but rotating and rescaling $S_{2,t}$ gives a continous path $A_t$ in $GL_2(\mathbb{R})$ such that $S_{1,t}$ and $A_t.S_ {2,t}$ satisfy that condition. Note that we necessary have $A_2.S_{2,2}=S_{2,2}$. Therefore $f^{-1}\bigl(S_{1,t},A_2.S_{2,t}\bigr)$ is a continuous path in $U$ that joins $S^1$ to $S^{2}$. So the subset $U$ is connected. Note that the  connectedness of $U$ clearly implies the connectedness of $U\cap Q_1(k_1^ {\alpha_1},\ldots,k_r^{\alpha_r}, -1^s)$.

The cases c) and d) are analogous and left to the reader.
\end{proof}

\begin{cor}
Let $Q(k_1,\ldots,k_r)$ be a stratum of quadratic differentials on $\mathbb{CP}^1$, and let $N\geq 1$.
If a connected component of $Q^N(k_1,\ldots,k_r)$ admit orbifoldic points, then the corresponding configuration is symmetric and the locus of orbifoldic points are unions of copies (or coverings) of open subset of configuration domains of other strata, which are manifolds.
\end{cor}

\begin{proof}
Recall that $S$ correspond to an orbifoldic point if and only if $S$ admits a nontrivial orientation preserving isometry. Now let $U$ be a connected component of $\mathcal{Q}^N(k_1^{\alpha_1},\ldots,k_r^{\alpha_r}, -1^s)$, $S\in U$ an orbifoldic point, and $\tau$ an isometry.

Suppose that $U$ correspond to the $a)$ case. Then $\tau$ must preserve the smallest saddle connection $\gamma_0$ of $S$. Either $\tau$ fixes the end points of $S$, either it interchanges them. In the first case, $\tau=Id$,  in the other case it is uniquely determined and is an involution that fixes the middle of $\gamma_0$. In that case the end points of $\gamma_0$ have the same order $k\geq 0$. Then $S/ \tau$ is a half-translation surface whose smallest saddle connection is of length $|\gamma_0|/2$ and joins a singularity of order $k\geq 0$ to a pole. Any other saddle connection in $S/ \tau$ is of length $l$ or $l/2$ for $l$ the length of a saddle connection (different from $\gamma_0$) on $S$. Therefore, $S/ \tau$ belongs to a configuration domain of $a)$ type in the corresponding stratum. The flat surface  $S/ \tau$ does not have a nontrivial orientation preserving isometry because $k\neq -1$.  Therefore the configuration domain that contains $S/\tau$ is a manifold. The involution $\tau$ induces an involution on the set of zeros of $S$ and the stratum and configuration domain corresponding to $S/\tau$ depends only on that involution. This induces a covering from the locus of orbifoldic points whose corresponding involution share the same combinatorial data to an open set of a manifold. 

If $U$ belongs to the $b)$ case, then similarly, a nontrivial isometric involution $\tau$ interchanges the two 1-holed sphere of the decomposition.  We have $A_1=A_2$ and $a_1=a_2>0$ (see notations of Theorem \ref{th:config}). The set of orbifoldic points is isomorphic to the configuration domain of $a)$ type with data $\{a_1,-1\}$ which is a manifold.  

 If $U$ belongs to the $c)$ case then similarly, $\tau$ interchanges the two 1-holed sphere of the decomposition.  We must have $A_1=A_2$ and $a_1=a_2>0$. 
The set of orbifoldic points is isomorphic to an open subset of a configuration domain of $d)$ type,  which is a manifold (see next).  

In the $d)$ case, any isometry $\tau$ fix the saddle connection $\gamma_1$ that separates the surface in a 1-holed sphere and a half-pillowcase, which are nonisometric. Hence they are fixed by  $\tau$. Now since $\tau$ is orientation preserving, it is easy to check that necessary, $\tau$ is trivial.


\end{proof}

Here we use Theorem \ref{cusp_elementaire} and the description of configurations  to show that any stratum of quadratic differentials on $\mathbb{CP}^1$ admits only one topological end.

\begin{prop}\label{bout}
Let $\mathcal{Q}_1(k_1,\ldots,k_r,-1^s)$ be any strata of quadratic differential on $\mathbb{CP}^1$. Then the subset $\mathcal{Q}_{1,\delta}(k_1,\ldots,k_r,-1^s)$ is connected for any $\delta>0$.
\end{prop}

\begin{proof}

Let $S\in \mathcal{Q}_{1,\delta} (k_1,\ldots,k_r,-1^s)$. We first describe a path  from $S$ to a simple configuration domain with corresponding singularities of orders $\{-1,k\}$. Then we show that all of these configuration domains are in the same connected component of $\mathcal{Q}_{1,\delta} (k_1,\ldots,k_r,-1^s)$.

Let $\gamma_1$ be a saddle connection of $S$ of length less than $\delta$. Up to the Teichmüller geodesic flow action, we can assume that $\gamma_1$ is of length less than $\delta^2$. Now let $P$ be a pole. There exists a saddle connection $\gamma_2$ of length less than $1$ starting from $P$, otherwise the $1-$neighborhood of $P$ would be an embedded half-disk of radius $1$  in the surface, and would be of area $\frac{\pi}{2}>1$. Then up to a slight deformation, we can assume that there are no other saddle connections parallel to $\gamma_1$ or $\gamma_2$ (except the ones that are \^homologous to $\gamma_1$ or $\gamma_2$). 
Now we contract $\gamma_2$ using the Teichmüller geodesic flow. This 
gives a path in $\mathcal{Q}_{1,\delta} (k_1,\ldots,k_r,-1^s)$, and we now can assume that $\gamma_2$ is of length smaller than $\delta$.

The other end of $\gamma_2$ is a singularity of order $k$. If $k\geq 0$, then from the list of configurations given in Theorem \ref{th:config}, the saddle connection $\gamma_2$ is simple.

We assume that $k=-1$, then the surface is a  $1-$holed sphere glued with a cylinder, one end of this cylinder is $\gamma_2$ (we have a half-pillowcase), and the other end of that cylinder is a closed saddle connection whose end point is a singularity $P^{\prime}$ of order $k^{\prime}>0$. We can assume, up to applying the Teichmüller geodesic flow, that $\gamma_2$ is of length at most $(1-c)\delta$, where $c$ is the area of the cylinder.
 Now we consider $\gamma_3$ to be the shortest path from $P$ to $P^{\prime}$. It is clear that $\gamma_3$ is a simple saddle connection. Now up to twisting and shrinking the cylinder, we can make this saddle connection as small as possible (see Figure \ref{un:seul:bout}). However, this transformation, is not area preserving and we must rescale the surfaces to keep area one surfaces . This rescalling increase the length of $\gamma_2$ by a factor which is at most $\frac{1}{1-c}$, and therefore the length of $\gamma_2$ is always smaller than $\delta$ during this last deformation, and the resulting surface is in a simple  configuration domain with corresponding singularities of orders $\{-1,k^{\prime}\}$.
 
\begin{figure}
\begin{center}
\input{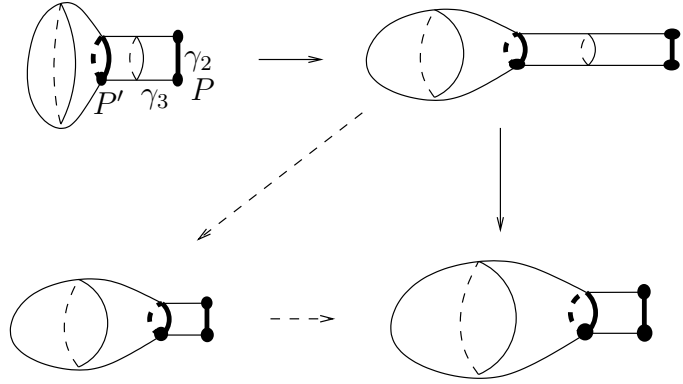}
\caption{Deformation of a surface in  $\mathcal{Q}_{1,\delta} (k_1,\ldots,k_r,-1^s)$.}
\label{un:seul:bout}
\end{center}
\end{figure}

Now let $(U_i)_{i=1,2}$ be simple configuration domains whose corresponding configurations are represented by simple paths that joins a pole to a singularity of order $k_i>0$ (here we assume that $r\geq 2$, the case $r\leq 1$ is trivial).  From Theorem \ref{cusp_elementaire}, for each $i=1,2$, the set $U_i \cap \mathcal{Q}_{1,\delta} (k_{1},\dots,k_{r},-1^{s})$ is connected. So, it is enough to find a path between two specific surfaces in $U_i$ that stays in $\mathcal{Q}_{1,\delta}(k_1,\ldots,k_r, -1^s)$. We start from a surface in $\mathcal{Q}(k_1-1,k_2-1,k_3,\ldots,k_r,-1^{s-2})$ and we successively break each singularity of order $k_i-1$ into two singularities of order $k_i$ and $-1$. We get a surface in $\mathcal{Q}_{1,\delta} (k_1,\ldots,k_r)$ with two arbitrary small saddle connections. We can assume that one of these short saddle connections is verticakl, and the other not. Then action on this surface by the Teichmüller geodesic flow easily give a path between $U_1$ and $U_2$ that keeps being in 
$\mathcal{Q}_{1,\delta}(k_1,\ldots,k_r, -1^s)$.
\end{proof}

\appendix
\section*{Appendix. A geometric criteria for \^homo\-logous saddle connections}
Here we give a proof of the following theorem:
\begin{NNths}[H. Masur, A. Zorich] 
Consider two distinct saddle connections $\gamma,\gamma^{\prime}$ on a half-translation surface. The following assertions are equivalent:
 \begin{itemize} 
 \item[a)]  The two saddle connections $\gamma$ and $\gamma^{\prime}$ are \^homologous.  
\item[b)] The ratio of their length is constant under any small deformation of the surface inside the ambient stratum.  
\item[c)]
They  have no interior intersection and one of the connected  component  of $S\backslash\{\gamma\cup \gamma^{\prime}\}$ has  trivial  linear holonomy.
\end{itemize} 
\end{NNths}

\begin{proof}
The the proofs of the statements $a \Leftrightarrow b$ and $c \Rightarrow b$ are easy and the arguments are same as in \cite{MZ}. We will write them for completeness. Our proof of $b \Rightarrow c$ is new and more geometric than the initial proof.

\begin{itemize}
\item
We first show that statement $a)$ is equivalent to statement $b)$.  We have defined $[\hat{\gamma}]$ and $[\hat{\gamma}^{\prime}]$ in $H_1^-(\widehat{S},\widehat{P},\mathbb{Z})$. They both are cycles associated to simple paths, therefore, they are primitive cycles of $H_1^-(\widehat{S},\widehat{P},\mathbb{Z})$. 

If $\gamma$ and $\gamma^{\prime}$ are \^homologous, then integrating $\omega$ along the cycles $[\hat{\gamma}]$ and $[\hat{\gamma^{\prime}}]$, we see that the ratio of their length belongs to $\{-1/2,1,2\}$, and this ratio is obviously constant under small deformations of the surface. Conversely, if they are not \^homolo\-gous, then $(\gamma,\gamma^{\prime})$ is a free family on  $H_1^-(\widehat{S},\widehat{P},\mathbb{C})$ (since they are primitive elements of $H_1^-(\widehat{S},\widehat{P},\mathbb{Z})$) and so $\int_{\hat{\gamma}} \omega$ and $\int_{\hat{\gamma^{\prime}}}\omega$ correspond to two independent coordinates in a neighborhood of $S$. Therefore the ratio of their length is not locally constant.

\item 
Now assume $c)$. We denote by $S^+$ a connected component of $S\backslash \{ \gamma,\gamma^{\prime}\}$ that has trivial holonomy.  Its boundary is a union of components homeomorphic to $\mathbb{S}^1$. The saddle connections have no interior intersections, so this boundary is a union of copies of $\gamma$ and $\gamma^{\prime}$ and it is easy to check that both $\gamma$ and $\gamma^{\prime}$ appears in that boundary. The flat structure on $S^+$ is defined by an Abelian differential $\omega$. Now we have  $\int_{\partial S^+} \omega =0$, which impose a relation on $|\gamma|$ and $|\gamma^{\prime}|$. This relation is preserved in a neighborhood of $S$, and therefore, the ratio is locally constant and belongs to $\{1/2,1,2\}$ depending on the number of copies of each saddle connections on the boundary os $S^+$.

\item
Now assume $b)$. We can assume that  the saddle connection $\sigma$ is vertical.  Then applying the Teichmüller geodesic flow $g_t$ on $S$, for some small $t$, induce a small deformation of $S$ . The hypothesis implies that the  saddle connection $\gamma^{\prime}$ is necessary vertical too, and so the two saddle connections are parallel and  hence have  no interior intersections. Let $S_1$ and $S_2$ the connected components of $S\backslash\{\gamma,\gamma^{\prime}\}$ that bounds $\gamma$ (we may have $S_1=S_2$), and assume that $S_1$ has nontrivial linear holonomy. That implies  there exists a simple broken line $\nu$ with nontrivial linear holonomy that starts and ends on the boundary of $S_1$ that correspond to $\gamma$. Now, we create an small hole by adding a parallelogram on the first segment of the path $\nu$. This create only one hole $\tau$ in the interior of $S_1$ because the other one is sent to the boundary (this procedure add the length of the hole to the length of the boundary). If we directly move the hole $\tau$ to the boundary, we obtain a flat surface isometric to  the initial surface $S_1$. But if we first transport $\tau$ along $\nu$, then this will change its orientation, and its length will be added again to the length of the boundary. So the resulting surface have a boundary component corresponding to $\gamma$ bigger than the initial surface $S_1$. The surgery did not affect the boundary corresponding to $\gamma^{\prime}$. Assume now that $S_2$ has also nontrivial holonomy, then performing the same surgery on $S_2$, and gluing back $S_1$ and $S_2$, this gives a slight deformation of $S$ that change the length of $\gamma$ and not the length of $\gamma^{\prime}$. This contradict the hypothesis $b)$.

\end{itemize}

\end{proof}

\end{document}